\documentstyle[amscd]{amsart}
\newsymbol\boxtimes 1202

\newcommand{\nc}{\newcommand}

\nc{\ad}{{{\bf{ad}}}}
\nc{\AJ}{{\operatorname{aj}}}
\nc{\Aut}{{\operatorname{Aut}}}
\nc{\Bls}{{{\cal B}ls}}
\nc{\Boxtimes}{{\fbox{$\times$}}}
\nc{\blt}{{\bullet}}
\nc{\bSt}{{\mbox{\bf{St}}}}
\nc{\card}{{\operatorname{card}}}
\nc{\Cch}{{\check{C}}}
\nc{\cd}{{\operatorname{cd}}}
\nc{\Ch}{{\operatorname{Ch}}}
\nc{\chara}{{\operatorname{char}}}
\nc{\CHom}{{\cal{H}om}}
\nc{\Coker}{{\operatorname{Coker}}}
\nc{\codim}{{\operatorname{codim}}}
\nc{\Cone}{{\operatorname{Cone}}}
\nc{\cSgn}{{\cal{S}gn}}
\nc{\depth}{{\operatorname{depth}}}
\nc{\dirlim}{{\underset{\rightarrow}{\operatorname{lim}}}}
\nc{\dotbox}{{\overset{\bullet}{\boxtimes}}}
\nc{\dotimes}{{\overset{\bullet}{\otimes}}}
\nc{\Ed}{{\operatorname{Edge}}}
\nc{\emp}{{\emptyset}}
\nc{\Ext}{{\operatorname{Ext}}}
\nc{\Fac}{{\cal{F}ac}}
\nc{\Fun}{{\operatorname{F}}}
\nc{\FS}{{\cal{FS}}}
\nc{\Hom}{{\operatorname{Hom}}}
\nc{\had}{{{\hat{\mbox{\bf{ad}}}}}}
\nc{\hgt}{{\operatorname{ht}}}
\nc{\Id}{{\operatorname{Id}}}
\nc{\id}{{\operatorname{id}}}
\nc{\Ima}{{\operatorname{Im}}}
\nc{\ind}{{\operatorname{ind}}}
\nc{\Ind}{{\operatorname{Ind}}}
\nc{\infi}{{\operatorname{inf}}}
\nc{\infh}{{\frac{\infty}{2}}}
\nc{\invlim}{{\underset{\leftarrow}{\operatorname{lim}}}}
\nc{\Jac}{{{\cal J}ac}}
\nc{\Ker}{{\operatorname{Ker}}}
\nc{\lcm}{{\operatorname{lcm}}}
\nc{\Locsys}{{{\cal L}ocsys}}
\nc{\Map}{{{\cal M}ap}}
\nc{\modul}{{\operatorname{mod}}}
\nc{\Mor}{{\operatorname{Mor}}}
\nc{\MS}{{\cal{MS}}}
\nc{\Ob}{{\operatorname{Ob}}}
\nc{\opp}{{\operatorname{opp}}}
\nc{\Or}{{{\cal O}r}}
\nc{\Ord}{{{\cal O}rd}}
\nc{\Part}{{{\cal P}art}}
\nc{\PGL}{{\operatorname{PGL}}}
\nc{\Pic}{{\operatorname{Pic}}}
\nc{\Rep}{{{\cal{R}}ep}}
\nc{\rk}{{\operatorname{rk}}}
\nc{\Sets}{{{\cal{S}}ets}}
\nc{\Sew}{{{\cal{S}}ew}}
\nc{\sgn}{{\operatorname{sgn}}}
\nc{\Sh}{{{\cal S}h}}
\nc{\Sign}{{{\cal S}ign}}
\nc{\Spe}{{\mbox{\bf{Sp}}}}
\nc{\supr}{{\operatorname{sup}}}
\nc{\Supp}{{\operatorname{Supp}}}
\nc{\supp}{{\operatorname{supp}}}
\nc{\Teich}{{{\cal{T}}eich}}
\nc{\tFS}{{\widetilde{\cal{FS}}}}
\nc{\Tor}{{\operatorname{Tor}}}
\nc{\totimes}{{\tilde{\otimes}}}
\nc{\tr}{{\operatorname{tr}}}
\nc{\tRep}{{\widetilde{{\cal R}ep}}}
\nc{\tTeich}{{\widetilde{{\cal T}eich}}}
\nc{\Vect}{{{\cal V}ect}}
\nc{\Ve}{{\operatorname{Vert}}}
\nc{\wt}{{\widetilde}}

\nc{\bo}{{{\bf{0}}}}
\nc{\One}{{{\bf{1}}}}
\nc{\one}{{{\bf{1}}}}

\nc{\BA}{{\Bbb A}}
\nc{\bA}{{\bf{A}}}
\nc{\bbA}{{\overline{\bf{A}}}}
\nc{\ba}{{{\bf{a}}}}
\nc{\baB}{{\overline{B}}}
\nc{\baeta}{{\bar{\eta}}}
\nc{\baJ}{{\bar{J}}}
\nc{\BB}{{\Bbb B}}
\nc{\bB}{{{\bf{B}}}}
\nc{\bc}{{{\bf{c}}}}
\nc{\bC}{{\overline{C}}}
\nc{\BC}{{\Bbb{C}}}
\nc{\bCC}{{\overline{\cal{C}}}}
\nc{\bCM}{{\overline{\cal{M}}}}
\nc{\bD}{{\bar{D}}}
\nc{\BD}{{\overline{D}}}
\nc{\bd}{{{\bf{d}}}}
\nc{\BE}{{\overline{E}}}
\nc{\BF}{{\Bbb{F}}}
\nc{\bF}{{{\bf{F}}}}
\nc{\bg}{{{\bf{g}}}}
\nc{\bG}{{{\bf{G}}}}
\nc{\BG}{{\Bbb G}}
\nc{\bGamma}{{\overline{\Gamma}}}
\nc{\bbH}{{     {{\bf{H}}}_a       }}
\nc{\bH}{{{\bf{H}}}}
\nc{\bi}{{{\bf{i}}}}
\nc{\bI}{{{\bf{I}}}}
\nc{\bj}{{{\bf{j}}}}
\nc{\bK}{{{\bf{K}}}}
\nc{\bL}{{{\bf{L}}}}
\nc{\BL}{{\Bbb{L}}}
\nc{\blambda}{{\bar{\lambda}}}
\nc{\bM}{{{\bf{M}}}}
\nc{\bmu}{{\vec{\mu}}}
\nc{\bN}{{{\bf{N}}}}
\nc{\BN}{{\Bbb{N}}}
\nc{\bnu}{{{\boldmath{${\nu}$}}}}
\nc{\bof}{{{\bf{f}}}}
\nc{\bp}{{{\bf{p}}}}
\nc{\BP}{{\Bbb P}}
\nc{\bP}{{{\bf{P}}}}
\nc{\BPO}{{\overset{\circ}{\BP}}}
\nc{\BQ}{{\Bbb Q}}
\nc{\bQ}{{\bf Q}}
\nc{\bq}{{{\bf{q}}}}
\nc{\BR}{{\Bbb{R}}}
\nc{\bR}{{{\bf{R}}}}
\nc{\br}{{{\bf{r}}}}
\nc{\breta}{{\bar{\eta}}}
\nc{\bs}{{{\bf{s}}}}
\nc{\bS}{{{\bf{S}}}}
\nc{\bt}{{{\bf{t}}}}
\nc{\bU}{{{\bf{U}}}}
\nc{\bV}{{{\bf{V}}}}
\nc{\bu}{{{\bf{u}}}}
\nc{\BUpsilon}{{\bar{\Upsilon}}}
\nc{\bw}{{{\bf{w}}}}
\nc{\bx}{{{\bf{x}}}}
\nc{\bX}{{{\bf{X}}}}
\nc{\BZ}{{\Bbb{Z}}}
\nc{\bz}{{{\bf{z}}}}
\nc{\bZ}{{{\bf{Z}}}}
\nc{\bzero}{{\boldmath{$0$}}}

\nc{\CA}{{\cal A}}
\nc{\CAD}{{\overset{\bullet}{\cal{A}}}}
\nc{\CAO}{{\overset{\circ}{\cal{A}}}}
\nc{\CB}{{\cal B}}
\nc{\CC}{{\cal C}}
\nc{\CalD}{{\cal D}}
\nc{\CE}{{\cal E}}
\nc{\CF}{{\cal F}}
\nc{\CG}{{\cal G}}
\nc{\CH}{{\cal H}}
\nc{\CI}{{\cal I}}
\nc{\CID}{{\overset{\bullet}{\cal{I}}}}
\nc{\CJ}{{\cal J}}
\nc{\CK}{{\cal K}}
\nc{\CL}{{\cal L}}
\nc{\CM}{{\cal M}}
\nc{\CN}{{\cal N}}
\nc{\CO}{{\cal O}}
\nc{\CP}{{\cal P}}
\nc{\CPO}{{\overset{\circ}{\cal{P}}}}
\nc{\CQ}{{\cal Q}}
\nc{\CR}{{\cal R}}
\nc{\CS}{{\cal S}}
\nc{\CT}{{\cal T}}
\nc{\CTD}{{\overset{\bullet}{\cal{T}}}}
\nc{\CTPO}{{\overset{\circ}{\cal{T}\cal{P}}}}
\nc{\CU}{{\cal{U}}}
\nc{\CV}{{\cal V}}
\nc{\CW}{{\cal W}}
\nc{\CX}{{\cal X}}
\nc{\CY}{{\cal Y}}
\nc{\CZ}{{\cal Z}}

\nc{\dCL}{{\overset{\bullet}{\cal{L}}}}
\nc{\dd}{{\operatorname{d}}}
\nc{\ddelta}{{\overset{\bullet}{\delta}}}
\nc{\dfu}{{\overset{\bullet}{\frak{u}}}}
\nc{\dlambda}{{\overset{\bullet}{\lambda}}}
\nc{\DO}{{\overset{\circ}{D}}}
\nc{\dpar}{{\partial}}
\nc{\dS}{{\overset{\bullet}{S}}}
\nc{\dT}{{\overset{\bullet}{T}}}

\nc{\fa}{{\frak{a}}}
\nc{\fA}{{\frak{A}}}
\nc{\fb}{{\frak{b}}}
\nc{\fB}{{\frak{B}}}
\nc{\fC}{{\frak{C}}}
\nc{\fD}{{\frak{D}}}
\nc{\fE}{{\frak{E}}}
\nc{\fF}{{\frak{F}}}
\nc{\ff}{{\frak{f}}}
\nc{\fg}{{\frak{g}}}
\nc{\fG}{{\frak{G}}}
\nc{\fH}{{\frak{H}}}
\nc{\fii}{{\frak{i}}}
\nc{\fj}{{\frak{j}}}
\nc{\fl}{{\frak{l}}}
\nc{\fL}{{\frak{L}}}
\nc{\fM}{{\frak{M}}}
\nc{\fN}{{\frak{N}}}
\nc{\fn}{{\frak{n}}}
\nc{\fp}{{\frak{p}}}
\nc{\fQ}{{\frak{Q}}}
\nc{\fu}{{\frak{u}}}
\nc{\fU}{{\frak{U}}}
\nc{\fV}{{\frak{V}}}
\nc{\fW}{{\frak{W}}}
\nc{\fZ}{{\frak{Z}}}

\nc{\hCH}{{\hat{\cal{H}}}}
\nc{\hCI}{{\hat{\cal{I}}}}
\nc{\hfC}{{\hat{\frak{C}}}}
\nc{\hfg}{{\hat{\frak{g}}}}
\nc{\hL}{{\hat{L}}}
\nc{\HO}{{\overset{\circ}{H}}}
\nc{\hpsi}{{\hat{\psi}}}
\nc{\hx}{{\hat{x}}}

\nc{\jo}{{\overset{\circ}{j}}}

\nc{\phid}{{\overset{\bullet}{\phi}}}

\nc{\tA}{{\tilde{A}}}
\nc{\ta}{{\tilde{a}}}
\nc{\tB}{{\tilde{B}}}
\nc{\tb}{{\tilde{b}}}
\nc{\tBP}{{\tilde{\BP}}}
\nc{\tC}{{\tilde{C}}}
\nc{\tc}{{\tilde{c}}}
\nc{\tCA}{{\tilde{\cal{A}}}}
\nc{\tCC}{{\tilde{\cal{C}}}}
\nc{\tCH}{{\tilde{\cal{H}}}}
\nc{\tCI}{{\tilde{\cal{I}}}}
\nc{\tCO}{{\tilde{\cal{O}}}}
\nc{\tCP}{{\tilde{\cal{P}}}}
\nc{\tCT}{{\tilde{\cal{T}}}}
\nc{\tD}{{\tilde{D}}}
\nc{\tDelta}{{\tilde{\Delta}}}
\nc{\tE}{{\tilde E}}
\nc{\tF}{{\tilde F}}
\nc{\tfD}{{\tilde{\frak{D}}}}
\nc{\tfF}{{\tilde{\frak{F}}}}
\nc{\tff}{{\tilde{\frak{f}}}}
\nc{\tfu}{{\tilde{\frak{u}}}}
\nc{\tJ}{{\tilde{J}}}
\nc{\tj}{{\tilde{j}}}
\nc{\tK}{{\tilde K}}
\nc{\tL}{{\tilde{L}}}
\nc{\tM}{{\tilde{M}}}
\nc{\tP}{{\tilde{P}}}
\nc{\tPhi}{{\tilde{\Phi}}}
\nc{\tpi}{\tilde{\pi}}
\nc{\TPO}{{\overset{\circ}{T\BP}}}
\nc{\tR}{{\tilde{R}}}
\nc{\tS}{{\tilde S}}
\nc{\tT}{{\tilde{T}}}
\nc{\ttau}{{\tilde{\tau}}}
\nc{\ttheta}{{\tilde{\theta}}}
\nc{\tU}{{\tilde{U}}}
\nc{\tUpsilon}{{\tilde{\Upsilon}}}
\nc{\tW}{{\tilde W}}
\nc{\ty}{{\tilde y}}
\nc{\tY}{{\tilde Y}}
\nc{\txi}{{\tilde{\xi}}}

\nc{\UD}{{\overset{\bullet}{U}}}
\nc{\UO}{{\overset{\circ}{U}}}

\nc{\vA}{{\vec{A}}}
\nc{\valpha}{{\vec{\alpha}}}
\nc{\vbeta}{{\vec{\beta}}}
\nc{\vc}{{\vec{c}}}
\nc{\vD}{{\vec{D}}}
\nc{\vd}{{\vec{d}}}
\nc{\vgamma}{{\vec{\gamma}}}
\nc{\vK}{{\vec{K}}}
\nc{\vlambda}{{\vec{\lambda}}}
\nc{\vmu}{{\vec{\mu}}}
\nc{\vnu}{{\vec{\nu}}}
\nc{\vo}{{\vec{0}}}
\nc{\vu}{{\vec{u}}}
\nc{\vx}{{\vec{x}}}
\nc{\vy}{\vec{y}}
\nc{\vzero}{\vec{0}}

\nc{\XO}{{\overset{\circ}{X}}}

\nc{\ya}{{\operatorname{aj}}}

\nc{\nen}{\newenvironment}
\nc{\ol}{\overline}
\nc{\ul}{\underline}
\nc{\ra}{\rightarrow}
\nc{\lra}{\longrightarrow}
\nc{\Lra}{\Longrightarrow}
\nc{\lla}{\longleftarrow}
\nc{\Llra}{\Longleftrightarrow}
\nc{\hra}{\hookrightarrow}
\nc{\iso}{\overset{\sim}{\lra}}
\nc{\rlh}{\rightleftharpoons}

\nc{\IC}{{\cal{IC}}}
\nc{\PS}{{\cal{PS}}}
\nc{\oCG}{{\overline{\cal G}}}
\nc{\oCQ}{{\overline{\cal Q}}}
\nc{\oCZ}{{\overline{\cal Z}}}
\nc{\dZ}{{\overset{\bullet}{\cal Z}}{}}
\nc{\dCS}{{\overset{\bullet}{\cal S}}{}}
\nc{\dbQ}{{\overset{\bullet}{\bf Q}}{}}
\nc{\ddZ}{{\ddot{\cal Z}}{}}
\nc{\oZ}{{\overset{\circ}{\cal Z}}{}}
\nc{\dP}{{\overset{\bullet}{\cal P}}{}}
\nc{\oP}{{\overset{\circ}{\cal P}}{}}
\nc{\oQ}{{\overset{\circ}{\cal Q}}{}}
\nc{\obp}{{\overset{\circ}{{\bf p}}}{}}
\nc{\tbj}{{\tilde{\bf j}}{}}
\nc{\tbp}{{\tilde{\bf p}}{}}
\nc{\tfC}{{\widetilde{\frak C}}{}}
\nc{\tfE}{{\widetilde{\frak E}}{}}
\nc{\tfj}{{\widetilde{\frak j}}{}}
\nc{\tbQ}{{\widetilde{\bf Q}}{}}
\nc{\hCQ}{{\widehat{\cal Q}}{}}
\nc{\tfQ}{{\widetilde{\frak Q}}{}}
\nc{\tfp}{{\widetilde{\frak p}}{}}
\nc{\ofQ}{{\overset{\circ}{{\frak Q}}}{}}
\nc{\tGQ}{{\widetilde{\cal{GQ}}}{}}
\nc{\tCS}{{\widetilde{\cal S}}{}}
\nc{\oGQ}{{\overset{\circ}{\cal{GQ}}}{}}
\nc{\ooGQ}{{\overset{\circ\circ}{\cal{GQ}}}{}}
\nc{\oGZ}{{\overset{\circ}{\cal{GZ}}}{}}
\nc{\tGZ}{{\widetilde{\cal{GZ}}}{}}
\nc{\ufM}{{\underline{\frak M}}}
\nc{\ufQ}{{\underline{\frak Q}}}
\nc{\Ue}{{U_\varepsilon}}
\nc{\Upe}{{\Upsilon_\varepsilon}}
\nc{\crho}{{\check{\rho}}}
\nc{\ctheta}{{\check{\theta}}}
\nc{\sZ}{{\sf Z}}

\nc{\Thm}[1]{Theorem~\ref{#1}}
\nc{\Prop}[1]{Proposition~\ref{#1}}
\nc{\Lem}[1]{Lemma~\ref{#1}}
\nc{\Cor}[1]{Corollary~\ref{#1}}
\nc{\Conj}[1]{Conjecture~\ref{#1}}
\nc{\Claim}[1]{Claim~\ref{#1}}
\nc{\Defn}[1]{Definition~\ref{#1}}
\nc{\Exa}[1]{Example~\ref{#1}}
\nc{\Rem}[1]{Remark~\ref{#1}}
\nc{\Note}[1]{Note~\ref{#1}}

\nen{thm}[1]{\label{#1}{\bf Theorem.\ } \em}{}
\nen{prop}[1]{\label{#1}{\bf Proposition.\ } \em}{}
\nen{lem}[1]{\label{#1}{\bf Lemma.\ } \em}{}
\nen{cor}[1]{\label{#1}{\bf Corollary.\ } \em}{}
\nen{conj}[1]{\label{#1}{\bf Conjecture.\ } \em}{}
\nen{claim}[1]{\label{#1}{\bf Claim.\ } \em}{}

\nen{defn}[1]{\label{#1}{\bf Definition.\ } }{}
\nen{exa}[1]{\label{#1}{\bf Example.\ } }{}

\nen{rem}[1]{\label{#1}{\em Remark.\ } }{}
\nen{note}[1]{\label{#1}{\em Note.\ } }{}
\nen{exer}[1]{\label{#1}{\em Exercise.\ } }{}

\begin{document}

%

\nc{\nn}{\newline}
\nc{\nnn}{\newpage}
\nc{\noi}{\noindent}
\nc{\nop}{\noindent {\bf Proof.} }
\nc{\sbr}{\smallpagebreak}
\nc{\mbr}{\medpagebreak}
\nc{\bbr}{\bigpagebreak}


\nc{\raa}{\longrightarrow}
\nc{\lala}{\leftarrow}
\nc{\laa}{\longleftarrow}
\nc{\lrax}{\leftrightarrow}     

\nc{\Ra}{\Rightarrow}         
\nc{\LRa}{\Rightarrow}        

\nc{\inj}{\hookrightarrow}    
\nc{\injj}{\hookleftarrow}    
\nc{\sur}{\twoheadrightarrow} 
\nc{\surr}{\twoheadleftarrow} 
\nc{\mm}{\mapsto}             
\nc{\va}{\uparrow}              


\nc{\bb}{\underset}           
\nc{\aax}{\overset}     

\nc{\bsx}{\backslash}           
\nc{\bss}{\backslash}           
\nc{\barr}{\overline}         
\nc{\sss}{\S}               

\nc{\sub}{\subseteq}          
\nc{\suppp}{\supseteq}          

\nc{\ti}{\tilde}              
\nc{\tii}{\widetilde}         
\nc{\ch}{\check}              

\nc{\tim}{\times}             
\nc{\btim}{\boxtimes}
\nc{\ten}{\otimes}            
\nc{\bten}{\boxtimes}         
\nc{\pl}{\oplus}              
\nc{\con}{ @>\cong>> }  
\nc{\conn}{     @<\cong<<  }    

\nc{\half}{ \frac{1}{2} }     

\nc{\ci}{\circ}               
\nc{\cdx}{\cdot}             
\nc{\cdd}{\cdot}             
\nc{\cddd}{\cdot\cdot\cdot}


\nc{\cupp}{\bigcup}             
\nc{\capp}{\bigcap}
\nc{\tenn}{\bigotimes}          
\nc{\pll}{\bigoplus}

\nc{\pii}{\prod}                
\nc{\ppii}{\bigprod}            

\nc{\cci}{\bigcoprod}
\nc{\wwe}{\bigwedge}            
\nc{\cce}{\bigcoprod}           

\nc{\pp}{\endproclaim}        
\nc{\hh}{\endheading}         

\nc{\df}{ \overset{ \text{def}}= }
\nc{\inv}{ {}^{-1}      }
\nc{\we}{\wedge}
\nc{\wee}{{     \overset{2}\wedge       }}
\nc{\ppp}{{ \Bbb P^1 }}            
\nc{\aaa}{{\Bbb A^1}}              
\nc{\qlb}{ \barr{\Q_l} }      
\nc{\ffq}{ {\F_q} }           
\nc{\Spec}{{\text{Spec}\ {}     }}
\nc{\aand}{ \ \text{and}\ }
\nc{\hk}{       \text{hyperk\"ahler}    }
\nc{\rank}{{\ \text{rank}\ }}
\nc{\timB} {{	\underset{B}\tim		}}
\nc{\timP}{{    \underset{P}\tim                }}
\nc{\timQ}{{    \underset{Q}\tim                }}
\nc{\ab}{       ^{ab}   }
\nc{\af}{       ^{aff}  }



\nc{\AAA}{\cal A}
\nc{\BBB}{\cal B}       
\nc{\DD}{\cal D}
\nc{\EE}{\cal E}
\nc{\FF}{\cal F}
\nc{\GG}{\cal G}
\nc{\HH}{\cal H}
\nc{\II}{\cal I}
\nc{\JJ}{\cal J}
\nc{\KK}{\cal K}
\nc{\LL}{\cal L}
\nc{\MM}{\cal M}
\nc{\NN}{\cal N}
\nc{\OO}{\cal O}
\nc{\PP}{\cal P}
\nc{\QQ}{\cal Q}
\nc{\RR}{\cal R}
\nc{\SSS}{\cal S}
\nc{\TT}{\cal T}
\nc{\UU}{\cal U}
\nc{\VV}{\cal V}
\nc{\ZZ}{\cal Z}
\nc{\XX}{\cal X}
\nc{\YY}{\cal Y}

\nc{\A}{\Bbb A }
\nc{\cs}{\Bbb C^*}
\nc{\ccs}{ \Bbb C^*}
\nc{\cc}{\Bbb C}
\nc{\f}{\Bbb F}
\nc{\g}{\Bbb G}
\nc{\h}{\Bbb H}
\nc{\I}{\Bbb I}
\nc{\J}{\Bbb J}
\nc{\K}{\Bbb K}
\nc{\M}{\Bbb M}
\nc{\N}{\Bbb N}
\nc{\p}{\Bbb P}
\nc{\Q}{\Bbb Q}
\nc{\R}{\Bbb R}
\nc{\s}{\Bbb S}
\nc{\T}{\Bbb T}
\nc{\U}{\Bbb U}
\nc{\V}{\Bbb V}
\nc{\Z}{\Bbb Z}
\nc{\X}{\Bbb X}
\nc{\Y}{\Bbb Y}

\nc{\fI}{\frak I}
\nc{\fJ}{\frak J}
\nc{\fK}{\frak K}
\nc{\fO}{\frak O}
\nc{\fP}{\frak P}
\nc{\fR}{\frak R}
\nc{\fS}{\frak S}
\nc{\fT}{\frak T}
\nc{\fX}{\frak X}
\nc{\fY}{\frak Y}
\nc{\fc}{\frak c}
\nc{\fd}{\frak d}
\nc{\fe}{\frak e}
\nc{\fh}{\frak h}
\nc{\fk}{\frak k}
\nc{\fm}{\frak m}
\nc{\fo}{\frak o}
\nc{\fq}{\frak q}
\nc{\fr}{\frak r}
\nc{\fs}{\frak s}
\nc{\ft}{\frak t}
\nc{\fv}{\frak v}
\nc{\fz}{\frak z}
\nc{\fx}{\frak x}
\nc{\fy}{\frak y}

\nc{\al}{\alpha }
\nc{\be}{\beta }
\nc{\ga}{\gamma }
\nc{\de}{\delta }
\nc{\del}{\partial }
\nc{\ep}{\varepsilon }
\nc{\vap}{\epsilon }

\nc{\ze}{\zeta }
\nc{\et}{\eta }
\nc{\vth}{\vartheta }

\nc{\io}{\iota }
\nc{\ka}{\kappa }
\nc{\la}{\lambda }
\nc{\vrho}{\varrho}
\nc{\si}{\sigma }
\nc{\ups}{\upsilon }
\nc{\vphi}{\varphi }
\nc{\om}{\omega }

\nc{\Ga}{\Gamma }
\nc{\De}{\Delta }
\nc{\Th}{\Theta }
\nc{\La}{\Lambda }
\nc{\Si}{\Sigma }
\nc{\Ups}{\Upsilon }
\nc{\Om}{\Omega }


\nc{\zp}{{\overset{\bullet}{\cal Z}}{}}
\nc{\zc}{{\overset{\circ}{\cal Z}}{}}
\nc{\qp}{{\overset{\bullet}{\cal Q }}{}}
\nc{\qc}{{\overset{\circ}{\cal Q}}{}}
\nc{\nii}{ ^{n\cdd i} }
\nc{\yy}{\infty}

\nc{\cz}{\cc[z]}
\nc{\czn}{{ \cc_{\le n}[z] }}
\nc{\pv}{{\p(V)}}
\nc{\qv}{Q(V)}
\nc{\qcv}{{\overset{\circ}Q(V)}}
\nc{\ppi}{\p(I)}
\nc{\pk}{\p(K)}
\nc{\zv}{\ZZ(V)}
\nc{\zcv}{\zc(V)}

\nc{\ii}{{i\in I}}
\nc{\all}{{ ^{(\alpha)} }}
\nc{\bee}{{ ^{(\beta)} }}
\nc{\gaa}{{ ^{(\gamma)} }}

\title[]{Semiinfinite flags. I. Case of global curve $\BP^1$.}
\author{Michael Finkelberg}
\address{Independent Moscow University, Bolshoj Vlasjevskij pereulok, dom 11,
Moscow 121002 Russia}
\email{fnklberg@@mccme.ru}
\author{Ivan Mirkovi\'c}
\address{Dept. of Mathematics and Statistics, University of Massachusetts
at Amherst, Amherst MA 01003-4515, USA}
\email{mirkovic@@math.umass.edu}
\thanks{M.F. is partially supported by the U.S. Civilian
Research and Development Foundation under Award No. RM1-265 and by
INTAS94-4720. I.M. is partially supported by NSF}
\maketitle

\section{Introduction}

\subsection{} We learnt of the {\em Semiinfinite Flag Space} from B.Feigin
and E.Frenkel in the late 80-s. Since then we tried to understand this
remarkable object. It appears that it was essentially constructed, but
under different disguises, by V.Drinfeld and G.Lusztig in the early 80-s.
Another recent discovery ({\em Beilinson-Drinfeld Grassmannian}) turned out
to conceal a new incarnation of Semiinfinite Flags.
We write down these and other results scattered in the folklore.

\subsection{}
Let $\bG$ be an almost simple simply-connected group with
a Cartan datum $(I,\cdot)$ and a simply-connected simple
root datum $(Y,X,\ldots)$ of finite type as in ~\cite{l}, ~2.2.
We fix a Borel subgroup $\bB\subset\bG$, with a Cartan subgroup
$\bH\subset\bB$, and the unipotent radical $\bN$.
B.Feigin and E.Frenkel define the Semiinfinite Flag Space
$\CZ$ as the quotient of 
$\bG((z))$ modulo the connected component of $\bB((z))$
(see ~\cite{ff}).
Then they study the category $\PS$ of perverse sheaves on $\CZ$ equivariant
with respect to the Iwahori subgroup $\bI\subset \bG[[z]]$.

In the first two chapters we are trying to make sense of this definition.
We encounter a number of versions of this space.
In order to give it a structure of an ind-scheme, we define the (local)
semiinfinite flag space as 
$
\widetilde{\bf Q}=
\bG((z))/\bH\bN((z))$ (see section 4).
The (global) semiinfinite space
attached to a smooth complete curve $C$ is the system of 
varieties $\CQ^\al$ of ``quasimaps'' from
$C$ to the flag variety of $\bG$
--- the Drinfeld compactifications of the degree $\al$ maps.
In the present work we restrict ourselves to the case $C=\BP^1$.

The main incarnation of the semiinfinite flag space  in this paper
is  a collection  $\CZ$ (for {\em zastava})
of (affine irreducible finite dimensional) algebraic varieties
$\CZ^\alpha_\chi\sub \CQ^\al$, together with certain closed embeddings and 
{\em factorizations}.
Our definition of $\CZ$ follows the scheme suggested by G.Lusztig in
~\cite{l2}, \S11: we approximate the ``closures'' of Iwahori orbits by their
intersections with the transversal orbits of the opposite Iwahori subgroup.
However, since the set-theoretic intersections of the above
``closures'' with the opposite Iwahori orbits can not be equipped with the
structure of algebraic varieties, we  postulate $\CZ^\alpha_\chi$
for the ``correct'' substitutes of such intersections.

Having got the collection of $\CZ^\alpha_\chi$ with factorizations, we
imitate the construction of ~\cite{fs} 
to define the category $\PS$ (for
{\em polubeskrajni snopovi}) 
of certain collections of perverse sheaves with $\BC$-coefficients
on $\CZ^\alpha_\chi$ with {\em factorization isomorphisms}. It is defined
in chapter 2; this category is the main character of the present work.

\subsection{} 
\label{quantum}
If $\bG$ is of type $A,D,E$ we set $d=1$; if $\bG$ is of type
$B,C,F$ we set $d=2$; if $\bG$ is of type $G_2$ we set $d=3$. Let $q$
be a root of unity of sufficiently large degree $\ell$ divisible by $2d$.
Let $\fu$ be the small (finite-dimensional) quantum group associated to $q$
and the root datum $(Y,X,\ldots)$ as in ~\cite{l}.
Let $\CC$ be the category of $X$-graded $\fu$-modules
as defined in ~\cite{ajs}. Let $\CC^0$ be the block of $\CC$ containing the
trivial $\fu$-module. B.Feigin and G.Lusztig conjectured (independently)
that the category $\CC^0$ is equivalent to $\PS$.

Let $\fU\supset\fu$ be the quantum group with divided powers associted to $q$
and the root datum $(Y,X,\ldots)$ as in ~\cite{l}. Let $\fC$ be the category
of $X$-graded finite dimensional $\fU$-modules, and let $\fC^0$ be the block
of $\fC$ containing the trivial $\fU$-module.
The works ~\cite{kl}, ~\cite{l4} and
~\cite{kt} establish an equivalence of $\fC^0$ and the category $\CP(\CG,\bI)$.
Here $\CG$ denotes the affine Grassmannian $\bG((z))/\bG[[z]]$, and
$\CP(\CG,\bI)$ stands for the category of perverse sheaves
on $\CG$ with finite-dimensional support constant along the orbits of $\bI$.

\subsection{}
\label{quantum res}
The chapter 3 is devoted to the construction of the {\em convolution} functor
$\bc_\CZ:\ \CP(\CG,\bI)\lra\PS$
which is the geometric counterpart
of the restriction functor from $\fC^0$ to $\CC^0$, as suggested by V.Ginzburg
(cf. ~\cite{gk} ~\S4).
One of the main results of this chapter is the Theorem ~\ref{Satake} which is
the sheaf-theoretic version of the classical Satake isomorphism. Recall that
one has a {\em Frobenius homomorphism} $\fU\lra U(\fg^L)$ (see ~\cite{l})
where $U(\fg^L)$ stands for the universal enveloping algebra of the Langlands
dual Lie algebra $\fg^L$. Thus the category of finite dimensional
$\bG^L$-modules is naturally embedded into $\fC$ (and in fact, into $\fC^0$).
On the geometric level this corresponds to the embedding $\CP(\CG,\bG[[z]])
\subset\CP(\CG,\bI)$. The Theorem ~\ref{Satake} gives a natural interpretation
(suggested by V.Ginzburg)
of the weight spaces of $\bG^L$-modules in terms of the composition
$$\bG^L-mod\simeq\CP(\CG,\bG[[z]])\subset\CP(\CG,\bI)\stackrel{\bc_\CZ}{\lra}
\PS.$$

\subsection{}
Let us also mention here the following conjecture which might be known to
specialists (characteristic $p$ analogue of conjecture in ~\ref{quantum}).
Let $\bG^L$ stand for the Langlands dual Lie group. Let
$p$ be a prime number bigger than the Coxeter number of $\fg^L$, and let
$\overline\BF_p$ be
the algebraic closure of finite field $\BF_p$.
Let $\fC_p$ be the category of algebraic $\bG^L(\overline\BF_p)$-modules,
and let $\fC^0_p$
be the block of $\fC_p$ containing the trivial module.
Let $\CC_p$ be the category of graded modules over the Frobenius kernel of
$\bG^L(\overline\BF_p)$, and let $\CC^0_p$
be the block of $\CC_p$ containing the trivial module (see ~\cite{ajs}).
Finally, let $\PS_p$ be the category of snops {\em with
coefficients in} $\overline\BF_p$, and let $\CP(\CG,\bI)_p$ be the category
of perverse sheaves on $\CG$ constant along $\bI$-orbits {\em with coefficients
in} $\overline\BF_p$. Then the categories $\CC^0_p$ and $\PS_p$
are equivalent, the categories $\fC^0_p$ and $\CP(\CG,\bI)_p$ are equivalent,
and under these equivalences the restriction functor $\fC^0_p\lra\CC^0_p$
corresponds to the convolution functor $\CP(\CG,\bI)_p\lra\PS_p$ (cf.
~\ref{quantum res}).
The equivalence $\CP(\CG,\bI)_p\iso\fC^0_p$ should be an extension of
the equivalence between $\CP(\CG,\bG[[z]])_p\subset\CP(\CG,\bI)_p$ and
the subcategory of $\fC^0_p$ formed by the
$\bG^L(\overline\BF_p)$-modules which factor through the Frobenius homomorphism
$Fr:\ \bG^L(\overline\BF_p)\lra\bG^L(\overline\BF_p)$. The latter equivalence
is the subject of forthcoming paper of K.Vilonen and the second author.

\subsection{}
The Zastava space $\CZ$ organizing all the ``transversal slices"
$\CZ^\alpha_\chi$ may seem cumbersome.
At any rate the existence of various models of the slices
$\CZ^\alpha_\chi$
(chapter 1), is undoubtedly beautiful by itself. Some of the
wonderful properties of $\CQ^\al$ and $\CZ^\alpha_\chi$
are demonstrated in ~\cite{ku}, ~\cite{fk}, ~\cite{fkm} in the case
$\bG=SL_n$. We expect all these properties to hold for the general $\bG$.

\subsection{}
To guide the patient reader through the
notation, let us list
the key points of this paper. The Theorem ~\ref{Z} identifies the
different models of $\CZ^\alpha_\chi$ (all essentially due to V.Drinfeld)
and states the factorization property. The exactness of the convolution
functor $\bc_\CZ:\ \CP(\CG,\bI)\lra\PS$ is proved in the
Theorem ~\ref{tough} and Corollary ~\ref{bunk}. The Theorem ~\ref{Satake}
computes the value of the convolution functor on $\bG[[z]]$-equivariant sheaves
modulo the parity vanishing conjecture ~\ref{parity}.

\subsection{}
In the next parts we plan to study
$D$-modules  on the local variety $\widetilde{\bf Q}$
(local construction of the category $\PS$, global sections as modules over
affine Lie algebra $\hat\fg$,
action of the affine Weyl group by Fourier transforms),
the relation of the local and global
varieties (local and global Whittaker sheaves, a version of the
convolution functor twisted by a character of $N((z))$),
and the sheaves on  Drinfeld compactifications of maps
into partial flag varieties.

\subsection{} The present work owes its very existence to V.Drinfeld. It
could not have appeared without the generous help of many people who shared
their ideas with the authors. 
Thus, the idea of {\em factorization} (section 9)
is due to V.Schechtman. 
A.Beilinson and V.Drinfeld taught us the
{\em Pl\"ucker} picture of the (Beilinson-Drinfeld) 
affine Grassmannian
(sections 6 and 10). 
G.Lusztig has computed the local singularities of the
Schubert strata closures in the spaces $\CZ^\alpha_\chi$ (unpublished,
cf ~\cite{l1}). 
B.Feigin and V.Ginzburg taught us their understanding of
the Semiinfinite Flags for many years (in fact, we learnt of Drinfeld's
Quasimaps' spaces from V.Ginzburg in the Summer 1995). 
A.Kuznetsov was always
ready to help us whenever we were stuck in the geometric problems (in fact,
for historical reasons, the section 3 has a lot in common with ~\cite{ku} \S1).
We have also benefited from the discussions with R.Bezrukavnikov and
M.Kapranov. 
Parts of this work were done while the  authors were enjoying
the hospitality and support of the University of Massachusetts at Amherst,
the Independent Moscow University and the Sveu\v{c}ili\v{s}te u Zagrebu.
It is a great pleasure to thank these institutions.

\section{Notations}

\subsection{}
\label{group}{\bf Group $\bG$ and its Weyl group $\CW_f$.}
We fix a Cartan datum $(I,\cdot)$ and a simply-connected simple
root datum $(Y,X,\ldots)$ of finite type as in ~\cite{l}, ~2.2.

Let $\bG$ be the corresponding simply-connected almost
simple Lie group with the Cartan subgroup $\bH$ and the Borel subgroup
$\bB\supset \bH$ corresponding to the set of simple roots $I\subset X$.
We will denote by $\CR^+\subset X$ the set of positive roots.
We will denote by $2\rho\in X$ the sum of all positive roots.

Let $\bB_+=\bB$ and  let $\bB_-\supset \bH$ be the opposite Borel subgroup.
Let $\bN$ (resp. $\bN_-$)
be the radical of $\bB$ (resp. $\bB_-$).
Let $\bbH=\bB/\bN=\bB_-/\bN_-$ be the abstract
Cartan group. The corresponding Lie algebras are denoted, respectively,
by $\fb,\fb_-,\fn,\fn_-,\fh$.

Let $\bX$ be the flag manifold $\bG/\bB$, and let $\bA=\bG/\bN$ be the
principal affine space. We have canonically
$H_2(\bX,\BZ)=Y;\ H^2(\bX,\BZ)=X$.

For $\nu\in X$ let $\bL_\nu$ denote the corresponding $\bG$-equivariant
line bundle on $\bX$.

Let $\CW_f$ be the Weyl group of $\bG$.
We have a canonical bijection $\bX^\bH=\CW_f$
such that the neutral element $e\in \CW_f=\bX^\bH\subset\bX$ forms a single
$\bB$-orbit.

We have a Schubert stratification of $\bX$ by $\bN$- (resp. $\bN_-$-)orbits:
$\bX=\sqcup_{w\in \CW_f}\bX_w$ (resp. $\bX=\sqcup_{w\in \CW_f}\bX^w)$ such that
for $w\in \CW_f=\bX^\bH\subset\bX$ we have $\bX^w\cap\bX_w=\{w\}$.

We denote by $\ol\bX_w$ (resp. $\ol\bX^w$) the Schubert variety --- the
closure of $\bX_w$ (resp. $\bX^w$). Note that $\ol\bX_w=\sqcup_{y\leq w}
\bX_y$ while $\ol\bX^w=\sqcup_{z\geq w}\bX^z$ where $\leq$ denotes the
standard Bruhat order on $\CW_f$.

Let $e\in \CW_f$ be the shortest element (neutral element),
let $w_0\in \CW_f$ be
the longest element, and let $s_i,\ i\in I$,
be the simple reflections in $\CW_f$.

\subsection{} 
\label{reps} {\bf Irreducible representations of $\bG$.}
We denote by $X^+$ the cone of positive weights (highest weights of finite
dimensional $\bG$-modules). The fundamental weights $\omega_i:\
\langle i,\omega_j\rangle=\delta_{ij}$ form the basis of $X^+$.

For $\lambda\in X^+$ we denote by $V_\lambda$ the finite dimensional
irreducible representation of $\bG$ with highest weight $\lambda$.

We denote by $V_\lambda^\vee$ the representation dual to $V_\lambda$;
the pairing: $V_\lambda^\vee\times V_\lambda\lra\BC$ is denoted by
$\langle,\rangle$.

For each $\lambda\in X^+$ we choose
a nonzero vector $y_\lambda\in V_\lambda^{\bN_-}$.
We also choose a nonzero vector $x_\lambda\in (V_\lambda^\vee)^\bN$ such that
$\langle x_\lambda,y_\lambda\rangle=1$.

\subsection{}
\label{config} {\bf Configurations of $I$-colored divisors.}
Let us fix $\alpha\in\BN[I]\subset Y,\ \alpha=\sum_{i\in I}a_ii$.
Given a curve $C$ we consider the configuration space $C^\alpha
\df\prod_{i\in I} C^{(a_i)}$ 
of colored
effective divisors of multidegree $\alpha$ (the set of colors is $I$).
The dimension of $C^\alpha$ is equal to the length
$|\alpha|=\sum_{i\in I}a_i$.


Multisubsets of a set $S$ are defined as elements of some symmetric power
$S^{(k)}$ and  we denote the image of $(s_1,...,s_k)\in S^k$ in $S^{(k)}$
by $\{\{s_1,...,s_k\}\}$.
We denote by $\fP(\alpha)$
the set of all partitions of $\alpha$, i.e multisubsets
$\Ga=
\{\{\ga_1,...,\ga_k\}\}$ of $\BN[I]$ with $\gamma_r\not=0$ and
$\sum_{r=1}^k \ga_i=\al$.

For $\Gamma\in\fP(\alpha)$ the corresponding stratum $C^\alpha_\Gamma$
is defined as follows. It is formed by configurations which can be
subdivided into $m$ groups of points, the $r$-th group containing $\gamma_r$
points; all the points in one group equal to each other, the different
groups being disjoint. For example, the main diagonal in $C^\alpha$
is the closed stratum given by partition $\alpha=\alpha$, while the complement
to all diagonals in $C^\alpha$ is the open stratum given by partition
$$
\alpha=\sum_{i\in I}(\underbrace{i_k+i_k+\ldots+i_k}_{a_k\operatorname{ times}})
$$
Evidently, $C^\alpha=\bigsqcup\limits_{\Gamma\in\fP(\alpha)}C^\alpha_\Gamma$.


\bigskip

\centerline{\bf  CHAPTER 1. The spaces $Q$ and $Z$}

\section{Quasimaps from a curve to a flag manifold}

\subsection{}
We fix a smooth projective curve $C$ and     $\alpha\in\BN[I]$.

\subsubsection{Definition} An algebraic map $f:\ C\lra\bX$ has degree $\alpha$
if the following equivalent conditions hold:

a) For the fundamental class $[C]\in H_2(C,\BZ)$ we have $f_*[C]=\alpha\in
Y=H_2(\bX,\BZ)$;

b) For any $\nu\in X$ the line bundle $f^*\bL_\nu$ on $C$ has degree
$\langle\alpha,\nu\rangle$.

\subsection{}
\label{maps}
The Pl\"ucker embedding of the flag manifold $\bX$ gives rise to the
following interpretation of algebraic maps of degree $\al$.

For any irreducible $V_\lambda$ we consider the trivial vector bundle
$\CV_\lambda=V_\lambda\otimes\CO$ over $C$.

For any $\bG$-morphism $\phi:\ V_\lambda\otimes V_\mu\lra V_\nu$ we denote
by the same letter the induced morphism $\phi:\ \CV_\lambda\otimes \CV_\mu
\lra \CV_\nu$.

Then
 a map  of degree $\al$ is a collection of {\em line subbundles}
$\fL_\lambda\subset\CV_\lambda,\ \lambda\in X^+$ such that:

a) $\deg\fL_\lambda=-\langle\alpha,\lambda\rangle$;

b) For any surjective
$\bG$-morphism $\phi:\ V_\lambda\otimes V_\mu\lra V_\nu$
such that $\nu=\lambda+\mu$ we have $\phi(\fL_\lambda\otimes\fL_\mu)=
\fL_\nu$;

c) For any $\bG$-morphism $\phi:\ V_\lambda\otimes V_\mu\lra V_\nu$
such that $\nu<\lambda+\mu$ we have $\phi(\fL_\lambda\otimes\fL_\mu)=0$.

Since the surjections
$V_\la\ten V_\mu\ra V_{\la+\mu}$ form
one $\cs$-orbit, systems $\LL_\la$ satisfying (b) are
determined by a choice of $\fL_{\om_i}$ for the  fundamental weights 
$\om_i,\ i\in I$.

If we replace the curve $C$ by a point, we get the Pl\"ucker description of
the flag variety $\bX$ as the set of collections of lines $L_\la\sub V_\la$
satisfying conditions of type (b) and (c).
Here, a Borel subgroup $B$ in $\bX$ corresponds to
a system  of lines $(L_\la,\ \la\in X^+)$ if the lines are the fixed
points of the unipotent radical $N$ of $B$, $L_\la=(V_\la)^N$,
or equivalently, if $N$ is the common stabilizer for all lines
$N=\bb{\la\in X^+}\cap G_{L_\la}$.

The space of degree $\alpha$ quasimaps from $C$ to $\bX$ will be denoted
by $\qc^\alpha$.

\subsection{Definition}
\label{quasimaps}
(V.Drinfeld)
The space
$\CQ^\al=\CQ^\al_C$ of {\em quasimaps} of degree $\alpha$
from $C$ to $\bX$ is the space of
collections of {\em invertible subsheaves}
$\fL_\lambda\subset\CV_\lambda,\ \lambda\in X^+$ such that:

a) $\deg\fL_\lambda=-\langle\alpha,\lambda\rangle$;

b) For any surjective 
$\bG$-morphism $\phi:\ V_\lambda\otimes V_\mu\lra V_\nu$
such that $\nu=\lambda+\mu$ we have $\phi(\fL_\lambda\otimes\fL_\mu)=
\fL_\nu$;

c) For any $\bG$-morphism $\phi:\ V_\lambda\otimes V_\mu\lra V_\nu$
such that $\nu<\lambda+\mu$ we have $\phi(\fL_\lambda\otimes\fL_\mu)=0$.

\subsubsection{Lemma}
a) The evident inclusion $\qc^\alpha\subset\CQ^\alpha$ is an open
embedding;

b) $\CQ^\alpha$ is a projective variety.

{\em Proof.} Obvious. $\Box$

\subsubsection{}
\label{praf}
Here is another version of the Definition, also due to V.Drinfeld.
The principal affine space $\bA=\bG/\bN$ is an $\bH_a$-torsor over $\bX$.
We consider its affine closure $\bbA$,
that is, the spectrum of the ring of functions on $\bA$.
Recall that $\bbA$ is the space of collections of vectors $v_\lambda\in
V_\lambda,\ \lambda\in X^+$, satisfying the following Pl\"ucker relations:

a) For any surjective
$\bG$-morphism $\phi:\ V_\lambda\otimes V_\mu\lra V_\nu$
such that $\nu=\lambda+\mu$, and $\phi(y_\lambda\otimes y_\mu)=y_\nu$,
we have $\phi(v_\lambda\otimes v_\mu)=v_\nu$;

b) For any $\bG$-morphism $\phi:\ V_\lambda\otimes V_\mu\lra V_\nu$
such that $\nu<\lambda+\mu$ we have $\phi(v_\lambda\otimes v_\mu)=0$.

The action of $\bH_a$ extends to $\bbA$ but it is not free anymore.
Consider the
quotient stack $\hat\bX=\bbA/\bH_a$. The flag variety $\bX$ is an
open substack in $\hat\bX$. A map
$\hat{\phi}:\ C\to\hat\bX$ is nothing else than an $\bH_a$-torsor
$\Phi$ over $C$ along
with an $\bH_a$-equivariant morphism $f:\ \Phi\to\bbA$. The degree of this map
is defined as follows.

Let $\lambda:\ \bH_a\to\BC^*$ be the character of $\bH_a$ corresponding
to a weight $\lambda\in X$. Let $\bH_\lambda\subset \bH_a$ be the kernel
of the morphism $\lambda$. Consider the induced $\BC^*$-torsor
$\Phi_\lambda=\Phi/\bH_\lambda$ over $C$. The map $\hat\phi$ has
degree $\alpha\in\BN[I]$ if
$$
\text{for any }\lambda\in X\quad\text{we have}\quad
\deg(\Phi_\lambda)=\langle\lambda,\alpha\rangle.
$$
{\bf Definition.}
The space $\CQ^\alpha$ is the space of maps $\hat{\phi}:\ C\to\hat\bX$
of degree $\alpha$ such that the generic point of $C$ maps into
$\bX\subset\hat\bX$.

The equivalence of the two versions of Definition follows by comparing their
Pl\"ucker descriptions.

\subsection{}
In this subsection we describe a stratification of $\CQ^\alpha$
according to the singularities of quasimaps.

\subsubsection{}
\label{sigma}
Given $\beta,\gamma\in\BN[I]$ such that $\beta+\gamma=\alpha$, we
define the proper map
$
\sigma_{\beta,\gamma}:\
\CQ^\beta	\times 	C^\gamma	\lra	\CQ^\alpha$.

Namely, let $f=(\fL_\lambda)_{\lambda\in X^+}\in\ \CQ^\beta$
be a quasimap of degree $\beta$;
and let
$
D=\sum_{i\in I}D_i
\cdd i
$
be an effective colored divisor of multidegree
$\gamma=\sum_{i\in I}d_ii$, that is, $\deg(D_i)=d_i$.
We define $\sigma_{\beta,\gamma}(f,D)\df f(-D)
\df	(
\fL_\lambda(-\langle D,\lambda\rangle )
)_{\lambda\in X^+}
\in\CQ^\alpha
,$
where we use the pairing $Div^I(C)\bb{\Z}\ten X\ra Div(C)$ given by
$\langle D,\lambda \rangle =
\sum_{i\in I}\langle i,\lambda\rangle \cdd D_i$.

\subsubsection{}
\label{strat M}
{\bf Theorem.} ${\displaystyle \ \CQ^\alpha=\bigsqcup_{0\leq\beta\leq\alpha}
\sigma_{\beta,\alpha-\beta}(\qc^\beta\times C^{\alpha-\beta})}$

{\em Proof.} Any invertible subsheaf $\fL_\la\sub \VV_\la$ lies in a unique
line subbundle $\ti\fL_\la\sub \VV_\la$ called the {\em normalization} of $\fL$.
So any quasimap $\fL$ defines a map $\ti\fL$ (called the {\em normalization}
of $\fL$) of degree $\be\le\al$
and an $I$-colored effective divisor $D$ (called the {\em defect} of $\fL$)
corresponding to the torsion sheaf $\ti\fL/\fL$, such that $\fL=\ti\fL(-D)$.
$\Box$

\subsubsection{Definition}
\label{domain}
Given a quasimap $f=(\fL_\lambda)_{\lambda\in X^+}
\in\ \CQ^\alpha$, its {\em domain of definition} $U(f)$
is the maximal Zariski open $U(f)\subset C$ such that for any $\lambda$
the invertible subsheaf $\fL_\lambda\subset\CV_\lambda$
restricted to $U(f)$ is actually a line subbundle.

\subsubsection{Corollary}
\label{big domain}
For a quasimap $f=(\fL_\lambda)_{\lambda\in X^+}\in\ \CQ^\alpha$ of degree
$\alpha$ the complement $C-U(f)$ of its domain of definition consists of
at most $|\alpha|$ points. $\Box$

\subsection{}
\label{C=line}
From now on, unless explicitly stated otherwise, $C=\BP^1$.

{\bf Proposition.} (V.Drinfeld) $\qc^\alpha$ is a smooth manifold of dimension
$2|\alpha|+\dim(\bX)$.

{\em Proof.} We have to check that at a map $f\in\qc^\alpha$ the first
cohomology $H^1(\BP^1,f^*\CT\bX)$ vanishes (where $\CT\bX$ stands for the
tangent bundle of $\bX$), and then the tangent space $\Theta_f\qc^\alpha$
equals $H^0(\BP^1,f^*\CT\bX)$.
As $\CT\bX$ is generated by the global sections, $f^*\CT\bX$ is
generated by global sections as well, hence
$H^1(\BP^1,f^*\CT\bX)=0$. To compute the dimension of $\Theta_f\qc^\alpha=
H^0(\BP^1,f^*\CT\bX)$ it remains to compute the Euler characteristic
$\chi(\BP^1,f^*\CT\bX)$. To this end we may replace $\CT\bX$ with its
associated graded bundle $\oplus_{\theta\in\CR^+}\bL_\theta$. Then
$$
\chi(\BP^1,f^*(\bigoplus_{\theta\in\CR^+}\bL_\theta))=\sum_{\theta\in\CR^+}
(\langle\alpha,\theta\rangle+1)=\langle\alpha,2\rho\rangle+\sharp\CR^+=
2|\alpha|+\dim\bX$$
$\Box$

\subsection{}
\label{CZ}
Now we are able to introduce our main character.
First we consider the open subspace $U^\alpha\subset\ \CQ^\alpha$ formed by
the quasimaps containing $\infty\in\BP^1$ in their domain of definition
(see ~\ref{domain}). Next we define the closed subspace $\CZ^\alpha\subset
U^\alpha$ formed by quasimaps with value at $\infty$ equal to $\bB_-\in\bX$:
$$\CZ^\alpha\df\ \{f\in U^\alpha | f(\infty)=\bB_-\}$$

We will see below that $\CZ^\alpha$ is an affine algebraic variety.

\subsubsection{}
\label{dimension}
It follows from Proposition ~\ref{C=line} that $\dim\CZ^\alpha=2|\alpha|$.

\section{Local Flag space}

In this section we define a version of $\CQ^\alpha$ where one replaces
the global curve $C$ by the formal neighbourhood of a point.

\subsection{}
\label{SS}
We set $\CO=\BC[[z]]\stackrel{p_n}{\lra}\CO_n=\BC[[z]]/z^n,\CK=\BC((z))$.

We define the scheme $\bbA(\CO)$ (of infinite type): its points are the
collections of vectors $v_\lambda\in V_\lambda\otimes\CO,\ \lambda\in X^+$,
satisfying the Pl\"ucker equations like in ~\ref{praf}.
It is a closed subscheme of $\prod_{i\in I}V_{\omega_i}\otimes\CO$.
We define the open subscheme $\bbA(\CO)_n\subset\bbA(\CO)$: it is formed
by the collections $(v_\lambda)_{\lambda\in X^+}$ such that $p_n(v_{\omega_i})
\not=0$ for all $i\in I$. Evidently, for $0\leq n\leq m$, one has
$\bbA(\CO)_n\subset\bbA(\CO)_m$.

We define the open subscheme $\CS\subset\bbA(\CO)$ as the union
$\bigcup_{n\geq0}\bbA(\CO)_n$. One has $\CS=\bA(\CO)$.

The scheme $\CS$ is equipped with the free action of $\bH_a:\
h(v_\lambda)_{\lambda\in X^+}=(\lambda(h)v_\lambda)_{\lambda\in X^+}$.
The quotient scheme $\bQ=\CS/\bH_a$ is a closed subscheme in
$\prod_{i\in I}\BP(V_{\omega_i}\otimes\CO)$. It is formed by the collections
of lines satisfying the Pl\"ucker equations. We denote the natural projection
$\CS\lra\bQ$ by $pr$.

\subsection{}
\label{SSeta}
For $\eta\in\BN[I]$ we define the closed subscheme $\CS^{-\eta}\subset\CS$
formed by the collections $(v_\lambda)_{\lambda\in X^+}$ such that
$v_\lambda=0\ \modul\ z^{\langle\eta,\lambda\rangle}$. We have the natural
isomorphism $\CS\iso\CS^{-\eta},\ (v_\lambda)_{\lambda\in X^+}\mapsto
(z^{\langle\eta,\lambda\rangle}v_\lambda)_{\lambda\in X^+}$.

Now we can extend the definition of $\CS^\chi$ to arbitrary $\chi\in Y$.
Namely, we define $\CS^\chi$ to be formed by the collections
$(v_\lambda\in V_\lambda\otimes\CK)_{\lambda\in X^+}$ such that
$(z^{\langle\chi,\lambda\rangle}v_\lambda)_{\lambda\in X^+}\in\CS$.
Evidently, $\CS^\chi\subset\CS^\eta$ iff $\chi\leq\eta$, and then the inclusion
is the closed embedding. The open subscheme $\CS^\eta-
\bigcup_{\chi<\eta}\CS^\chi\subset\CS^\eta$
will be denoted by $\dCS^\eta\subset\CS^\eta$.
The ind-scheme $\bigcup_{\eta\in Y}\CS^\eta$ will be denoted by $\tCS$.
The ind-scheme $\tCS$ is equipped with the natural action of the proalgebraic
group $\bG(\CO)$ (coming from the action on $\prod_{i\in I}V_{\omega_i}\otimes
\CK$), and the orbits are exactly $\dCS^\eta,\ \eta\in Y$.

\subsection{}
\label{QQ}
All the above (ind-)schemes are equipped with the free action
of $\bH_a$, and taking quotients we obtain the schemes
$\bQ^\eta=\CS^\eta/\bH_a,\ \eta\in Y$. They are all closed subschemes of
the ind-scheme $\prod_{i\in I}\BP(V_{\omega_i}\otimes\CK)$.
We have $\bQ^\chi\subset\bQ^\eta$ iff $\chi\leq\eta$, and then the inclusion
is the closed embedding. The ind-scheme $\tbQ=\tCS/\bH_a$ is the union
$\tbQ=\bigcup_{\eta\in Y}\bQ^\eta$.
The ind-scheme $\tbQ$ is equipped with the natural action of the proalgebraic
group $\bG(\CO)$ (coming from the action on
$\prod_{i\in I}\BP(V_{\omega_i}\otimes\CK)$),
and the orbits are exactly $\dbQ^\eta=\dCS^\eta/\bH_a,\ \eta\in Y$.

\subsection{}
We consider $C=\BP^1$ with two marked points $0,\infty\in C$. We choose a
coordinate $z$ on $C$ such that $z(0)=0,z(\infty)=\infty$.

\subsubsection{}
For $\alpha\in\BN[I]$ we define the space $\hCQ^\alpha\stackrel{pr}{\lra}
\CQ^\alpha$ formed by the collections $(v_\lambda\in\fL_\lambda\subset
\CV_\lambda)_{\lambda\in X^+}$ such that

a) $(\fL_\lambda\subset\CV_\lambda)_{\lambda\in X^+}\in\CQ^\alpha$;

b) $v_\lambda$ is a regular nonvanishing section of $\fL_\lambda$ on
$\BA^1=\BP^1-\infty$;

c) $(v_\lambda)_{\lambda\in X^+}$ satisfy the Pl\"ucker equations like in
~\ref{praf}.

It is easy to see that $\hCQ^\alpha\stackrel{pr}{\lra}\CQ^\alpha$ is a
$\bH_a$-torsor: $h(v_\lambda,\fL_\lambda)=(\lambda(h)v_\lambda,\fL_\lambda)$.

\subsubsection{}
\label{m}
Taking a formal expansion at $0\in C$ we obtain the closed embedding
$\fs_\alpha:\ \hCQ^\alpha\hookrightarrow\CS$. Evidently, $\fs_\alpha$ is
compatible with the $\bH_a$-action, so it descends to the same named
closed embedding $\fs_\alpha:\ \CQ^\alpha\hookrightarrow\bQ$.

\subsubsection{Lemma}
\label{codime}
Let $\beta\in\BN[I]$. Then $\codim_\bQ\bQ^{-\beta}\geq2|\beta|$.

{\em Proof.} Choose $\alpha\geq\beta$, and consider
the closed embedding $\fs_\alpha:\ \CQ^\alpha\hookrightarrow\bQ$.
Then $\fs_\alpha^{-1}(\bQ^{-\beta})=\CQ^{\alpha-\beta}$ embedded into
$\CQ^\alpha$ as follows: $(\fL_\lambda\subset\CV_\lambda)_{\lambda\in X^+}
\mapsto(\fL_\lambda(-\langle\beta,\lambda\rangle0)
\subset\CV_\lambda)_{\lambda\in X^+}$.
Now $\codim_\bQ\bQ^{-\beta}\geq\codim_{\CQ^\alpha}\CQ^{\alpha-\beta}=2|\beta|$.
$\Box$


\section{Pl\"ucker sections}

In this section we describe another model of the space $\CZ^\alpha$ introduced
in ~\ref{CZ}.

\subsection{}
\label{polynom}
We fix a coordinate $z$ on the affine line $\BA^1=\BP^1-\infty$.
We will also view  the configuration space
$\BA^\alpha\df(\BA^1)^\alpha$ (see ~\ref{config})
as the space of collections of unitary
polynomials $(Q_\lambda)_{\lambda\in X^+}$ in $z$, such that
(a) $\deg(Q_\lambda)=\langle\alpha,\lambda\rangle$, and
(b) $Q_{\lambda+\mu}=Q_\lambda Q_\mu$.

\subsection{}
\label{fZ}
Recall the notations of ~\ref{reps}.
For each $\lambda\in X^+$ we will use the decomposition
$V_\lambda=\BC y_\la\oplus(\Ker x_\la)=
(V_\lambda)^\bN\oplus\fn_- V_\lambda$,
compatible with the action of $\fh=\fb_-\cap\fb$, i.e., with the weight
decomposition.
For a section $v_\lambda\in\Ga(\BA^1,\CV_\lambda)=
V_\la\ten\BC[z]\df V_\la[z]$, we will use a polynomial
$Q_\la\df\ \langle x_\lambda,v_\lambda\rangle\in\BC[z]$,
to write down the decomposition
$v_\la=Q_\la\cdd y_\la\pl {v''}_\la\in
\BC[z]\cdd y_\la\pl(\Ker x_\la)[z] = V_\la[z]$.

{\bf Definition.} (V.Drinfeld) The space $\sZ^\alpha$ of
{\em Pl\"ucker sections} of degree $\alpha$ is the space of collections
of sections $v_\lambda\in\Ga(\BA^1,\CV_\lambda)=
V_\la\ten\BC[z]\df V_\la[z]
,\ \lambda\in X^+$;
such that for $v_\la=Q_\la\cdd y_\la\pl {v''}_\la\in
\BC[z]\cdd y_\la\pl(\Ker x_\la)[z]$, one has

a) Polynomial $Q_\lambda$ is unitary
of degree $\langle\alpha,\lambda\rangle$;

b) Component ${v''}_\la$ of $v_\la$ in $(\Ker x_\la)[z]$ has degree strictly
less than $\langle\alpha,\lambda\rangle$;

c) For any $\bG$-morphism $\phi:\ V_\lambda\otimes V_\mu\lra V_\nu$ such that
$\nu=\lambda+\mu$ and $\phi^\vee(x_\nu)=x_\lambda\otimes x_\mu$ we have
$\phi(v_\lambda\otimes v_\mu)=v_\nu$;

d) For any $\bG$-morphism $\phi:\ V_\lambda\otimes V_\mu\lra V_\nu$ such that
$\nu<\lambda+\mu$ we have $\phi(v_\lambda\otimes v_\mu)=0$.

\subsubsection{}
\label{affine}
Collections
$(v_\lambda)_{\lambda\in X^+}$ that satisfy (c),
are determined by a choice of $v_{\omega_i},\ i\in I$.
Hence $\sZ^\alpha$ is an affine algebraic variety.

\subsubsection{}
\label{pi}
Due to the properties a),c) above, the collection of polynomials $Q_\lambda$
defined in a) satisfies the conditions of  ~\ref{polynom}. Hence we
have the map $$\pi_\alpha:\ \sZ^\alpha\lra\BA^\alpha$$

\section{Beilinson-Drinfeld Grassmannian}

In this section we describe yet another model of the space $\CZ^\alpha$
introduced in ~\ref{CZ}.

\subsection{}
Let $C$ be an arbitrary smooth projective curve;
let $\CT$ be a left
$\bG$-torsor
over $C$, and let $\tau$ be a section of $\CT$ defined over a Zariski open
subset $U\subset C$, i.e., a trivialization of $\CT$ over $U$.
We will define a $\bB$- (resp. $\bB_-$-) type
$d(\tau)$ (resp. $d_-(\tau)$): a measure of singularity of $\tau$ at $C-U$.

\subsubsection{}
\label{type}
Section $\tau$ defines a $\bB$-subtorsor $\bB\cdd\tau\sub\CT$. This
reduction of $\CT$ to $\bB$ over $U$ is the same as a section of
$\bB\bss\CT$ over $U$. Since $\bG/\bB$ is
proper, this reduction (i.e. section),
extends uniquely to the whole $C$. Thus we obtain
a $\bB$-subtorsor
$\barr{\bB\cdd\tau}\sub\CT$ (the closure of $\bB\cdd\tau\sub\CT|U$ in $\CT$),
equipped with a section 
$\tau$ defined over $U$.

Using the projection $\bB\lra\bbH$ we can induce
$\barr{\bB\cdd\tau}$ to
a torsor over $C$ for the abstract  Cartan group $\bbH\cong \bB/\bN$ of $\bG$;
namely, $\CT_{\tau,\bB}\df \bN\bss \barr{\bB\cdd\tau}$, equipped with a section
$\tau_\bB$ defined over $U$.

The choice of simple coroots (cocharacters of $\bbH$) $I\subset Y$ identifies
$\bbH$ with $(\BC^*)^I$. Thus the section $\tau_\bB$
 of $\CT_{\tau,\bB}$ produces
an $I$-colored divisor $d(\tau)$ 
supported at $C-U$. We will call $d(\tau)$
the $\bB$-$type$ of $\tau$.

Replacing $\bB$ by $\bB_-$ in the above construction we define the
$\bB_-$-$type$
$d_-(\tau)$.

\subsection{}
\label{bZ}
Recall that A.Beilinson and V.Drinfeld have introduced the {\em relative
Grassmannian} $\CG_C^{(n)}$ over $C^n$ for any $n\in\BN$ (see
~\cite{todisappear}): its fiber $p_n^{-1}(x_1,\ldots,x_n)$
over an $n$-tuple $(x_1,\ldots,x_n)\in C^n$
is the space of isomorphism classes of $\bG$-torsors $\CT$ equipped with a
section $\tau$ defined over $C-\{x_1,\ldots,x_n\}$.

We will consider a certain finite-dimensional subspace of 
a partialy symmetrized version of the
relative Grassmannian.

{\bf Definition.}   (A.Beilinson and V.Drinfeld)
$\bZ^\alpha$ is the space of isomorphism classes of the following data:

a) an $I$-colored effective divisor $D\in\BA^\alpha$;

b) $\bG$-torsor $\CT$ over $\BP^1$ equipped with a section $\tau$ defined over
$\BP^1-supp(D)$ such that:

i) $\bB$-type $d(\tau)=0$;

ii) $\bB_-$-type $d_-(\tau)$ is a negative divisor (opposite to effective)
such that $d_-(\tau)+D$ is effective.

\subsubsection{}
\label{bZU}
By the definition, the space $\bZ^\alpha$ is equipped with a projection
$p_\alpha$ to $\BA^\alpha:\ (D,\CT,\tau)\mapsto D$. For a subset $U\subset\BA^1$
we will denote by $\bZ_U^\alpha$ the preimage $p_\alpha^{-1}(U)$.

\subsubsection{}
The reader may find another realization of $\bZ^\alpha$ in ~\ref{PBD} below.

\subsection{}
In this subsection we will formulate the crucial {\em factorization} property
of $\bZ^\alpha$.

\subsubsection{}
Recall the following property of the Beilinson-Drinfeld relative Grassmannian
$\CG_C^{(n)}\overset{p_n}{\lra}C^n$ (see ~\cite{todisappear}).
Suppose an $n$-tuple $(x_1,\ldots,x_n)\in C^n$ is represented as a union
of an $m$-tuple $(y_1,\ldots,y_m)\in C^m$ and
a $k$-tuple $(z_1,\ldots,z_k)\in C^k,\ k+m=n$, such that all the points of the
$m$-tuple are disjoint from all the points of the $k$-tuple.
Then $p_n^{-1}(x_1,\ldots,x_n)$ is canonically isomorphic to the product
$p_m^{-1}(y_1,\ldots,y_m)\times p_k^{-1}(z_1,\ldots,z_k)$

\subsubsection{}
\label{factorization}
Suppose we are given a decomposition $\alpha=\beta+\gamma,\ \beta,\gamma
\in\BN[I]$ and two disjoint subsets
$U,\Upsilon\subset\BA^1$.
Then $U^\beta\times\Upsilon^\gamma$ lies in $\BA^\alpha$, and we will denote
the preimage $p_\alpha^{-1}(U^\beta\times\Upsilon^\gamma)$ in $\bZ^\alpha$ by
$\bZ^{\beta,\gamma}_{U,\Upsilon}=
\bZ^\al|_{(U^\beta\times\Upsilon^\gamma)}$ (cf. ~\ref{bZU}).

The above property of relative Grassmannian immediately implies the
following

{\bf Factorization property.} There is a canonical factorization isomorphism
$\bZ^{\beta,\gamma}_{U,\Upsilon}\cong\bZ^\beta_U
\times\bZ^\gamma_\Upsilon$, i.e.,
$$
\bZ^\al |_{(U^\beta\times\Upsilon^\gamma)}\cong
\bZ^\be|_{U^\beta}
\times
\bZ^\ga|_{\Upsilon^\gamma}
.$$

\subsection{Remark} Let us describe the fibers of $p_\alpha$
in terms of the normal slices to the semiinfinite Schubert cells
in the loop Grassmannian.

\subsubsection{}
\label{Iwasawa}
Let $\CG$ be the usual affine Grassmannian $\bG((z))/\bG[[z]]$. It is
naturally identified with the fiber of $\CG^{(1)}_{\BP^1}$ over the point
$0\in\BP^1$. Due to the
Iwasawa decomposition in p-adic groups,
there is a natural bijection
between $Y$ and the set of orbits of the group
$\bN((z))$ (resp. $\bN_-((z))$) in $\CG$; for $\gamma\in Y$ we will denote the
corresponding orbit by $S_\gamma$ (resp. $T_\gamma$). We will denote by
$\ol{T}_\gamma$ the ``closure'' of $T_\gamma$, that is, the union
$\cup_{\beta\geq\gamma}T_\gamma$.

It is proved in ~\cite{mv} that the intersection
$\ol{T}_\gamma\cap S_\beta$ is not empty iff $\gamma\leq\beta$.
Then it is
an affine algebraic variety, a kind of a normal slice to
$T_\be$ in $\barr T_\ga$.
Let us call it $TS_{\gamma,\beta}
\df \ol{T}_\gamma\cap S_\beta$ for short.
If rank$(\bG)>1$ then $TS_{\gamma,\beta}=
\ol{T}_\gamma\cap S_\beta$ is not necessarily irreducible.
But it is always equidimensional of dimension $|\beta-\gamma|$. There is
a natural bijection between the set of irreducible components of
$TS_{\gamma,\beta}=\ol{T}_\gamma\cap S_\beta$
and the canonical basis of $U^+_{\beta-\gamma}$ (the
weight $\beta-\gamma$ component of the quantum universal enveloping algebra
of $\fn$) (see ~\cite{l} for the definition of canonical basis of $U^+$).

\subsubsection{}
Recall the diagonal stratification of $\BA^\alpha$ defined in ~\ref{config}
and the map $p_\al:\bZ^\al\ra\BA^\alpha$.
We consider a partition $\Gamma:\ \alpha=\sum_{k=1}^m\gamma_k$ and a divisor
$D$ in the stratum $\BA_\Gamma^\alpha$. The interested reader will check
readily the following

{\em Claim.} $p_\alpha^{-1}(D)$ is isomorphic to the product
$\prod_{k=1}^m TS_{-\gamma_k,0}
=\prod_{k=1}^m \ol{T}_{-\ga_k}\cap S_0\cong
\prod_{k=1}^m \ol{T}_0\cap S_{\ga_k}$.

In particular, the fiber over a point in the closed stratum is isomorphic to
$TS_{-\alpha,0}= \ol{T}_{-\al}\cap S_0 \cong  \ol{T}_0\cap S_\al $,
while the fiber over a generic point is isomorphic to
the product of affine lines
$TS_{-i,0}\cong \ol{T}_0\cap S_{-i}\cong\BA^1$, that is,
the affine space $\BA^{|\alpha|}$.

\subsubsection{Corollary}
\label{irred}
$\bZ^\alpha$ is irreducible.

\section{Equivalence of the three constructions}

\subsection{}
In this subsection we construct an isomorphism
$\varpi:\ \CZ^\alpha\iso\sZ^\alpha$, i.e., from the subsheaves
$\fL_\la\sub\CV_\la$ we construct the sections $v_\la\in\Ga(\BA^1,\CV_\la)$.

\subsubsection{}
\label{tuda}
Let $f\in\CZ^\alpha$ be a quasimap given by a collection $(\fL_\lambda
\subset\CV_\lambda=V_\lambda\otimes\CO_{\BP^1})_{\lambda\in X^+}$.
Since $\fL_\lambda|_{\BA^1}$ is trivial, it has a unique up to proportionality
section $v_\lambda$ generating it over $\BA^1$.

We claim that the pairing
$\langle x_\lambda,v_\lambda\rangle$ does not vanish identically. In effect,
since $\deg(f)=\alpha$, the meromorphic section
$\frac{v_\lambda}{z^{\langle\alpha,\lambda\rangle}}$
of $\CV_\lambda$ is regular nonvanishing
at $\infty\in\BP^1$. Moreover, since $f(\infty)=\bB_-$, we have
$\frac{v_\lambda}{z^{\langle\alpha,\lambda\rangle}}(\infty)\in
V_\lambda^{\bN_-}$.
Thus, $\langle x_\lambda,\frac{v_\lambda}{z^{\langle\alpha,\lambda\rangle}}
\rangle(\infty)\not=0$.

Now we can normalize $v_\lambda$ (so far defined up to a multiplication by
a constant) by the condition that $\langle x_\lambda,v_\lambda\rangle$ is a
unitary polynomial. Let us denote this polynomial by $Q_\lambda$. It has degree
$d_\lambda\leq\langle\alpha,\lambda\rangle$ since $\deg(f)=\alpha$.
Since
$\frac{v_\lambda}{z^{\langle\alpha,\lambda\rangle}}(\infty)\in
V_\lambda^{\bN_-}$,
we see that $\deg\langle e,v_\lambda\rangle<d_\lambda$ for any $e\perp y_\la$.
Moreover, since $\deg(f)=\alpha$ we must then have
$d_\lambda=\langle\alpha,\lambda\rangle$.

Thus we have checked that the collection $(v_\lambda)_{\lambda\in X^+}$
satisfies the conditions a),b) of the Definition ~\ref{fZ}. The conditions
c),d) of {\em loc. cit.} follow from the conditions b),c) of the Definition
~\ref{quasimaps}. In other words, we have defined the Pl\"ucker section
$$\varpi(f)\df\ (v_\lambda)_{\lambda\in X^+}\in\sZ^\alpha$$

\subsubsection{}
\label{obratno}
Here is the inverse construction. Given a Pl\"ucker section
$(v_\lambda)_{\lambda\in X^+}\in\sZ^\alpha$ we define the corresponding
quasimap $f=(\fL_\lambda)_{\lambda\in X^+}\in\CZ^\alpha$ as follows.

We can view $v_\lambda$ as a regular section of
$\CV_\lambda(\langle\alpha,\lambda\rangle\infty)$ over the whole $\BP^1$.
It generates an invertible subsheaf $\fL_\lambda'\subset
\CV_\lambda(\langle\alpha,\lambda\rangle\infty)$. We define
$$\fL_\lambda\df\ \fL_\lambda'(-\langle\alpha,\lambda\rangle\infty)
\subset\CV_\lambda$$

\subsubsection{}
It is immediate to see that the above constructions are inverse to each other,
so that $\varpi:\ \CZ^\alpha\lra\sZ^\alpha$ is an isomorphism.

\subsubsection{Remark} Note that the definition of the space $\CZ^\alpha$
depends only on the choice of Borel subgroup $\bB_-\subset \bG$, while the
definition of $\sZ^\alpha$ depends also on the choice of the opposite Borel
subgroup $\bB\subset \bG$ or, equivalently, on the choice of the Cartan subgroup
$\bH\subset \bB_-$.

We want to stress that the projection $\pi_\alpha:\ \CZ^\alpha=\sZ^\alpha\lra
\BA^\alpha$ {\em does depend} on the choice of $\bB$. Let us describe
$\pi_\alpha\varpi(f)$ for a genuine map (as opposed to quasimap)
$f\in\CZ^\alpha$. To this end recall (see ~\ref{group}) that the 
$\bB$-invariant
Schubert
varieties $\ol\bX_{s_iw_0},\ i\in I$, are divisors in $\bX$. Their formal
sum may be viewed as an $I$-colored divisor $\fD$ in $\bX$. Then $f^*\fD$ is
a well defined $I$-colored divisor on $\BP^1$ since $f(\BP^1)\not\subset\fD$
since $f(\infty)=\bB_-\in\bX_{w_0}$. For the same reason the point $\infty$
does not lie in $f^*\fD$, so $f^*\fD$ is really a divisor in $\BA^1$.
It is easy to see that $f^*\fD\in\BA^\alpha$ and $f^*\fD=\pi_\alpha\varpi(f)$.

\subsection{}
In this subsection we construct an isomorphism
$\xi:\ \sZ^\alpha\iso\bZ^\alpha$, so from a system of sections
$v_\la$ we construct a $\bG$-torsor $\CT$ with a section $\tau$
and an $I$-colored divisor $D$.

\subsubsection{Lemma}({\em The Pl\"ucker picture of $\bG$.})
\label{drinf}
The map
$\psi:\ g\mapsto (gx_\lambda,gy_\lambda)_{\lambda\in X^+}$ is a bijection
between $\bG$ and the space of collections $\{(u_\lambda\in V_\lambda^\vee,
\upsilon_\lambda\in V_\lambda)_{\lambda\in X^+})\}$ satisfying the following
conditions:

a) For any $\bG$-morphism $\phi:\ V_\lambda\otimes V_\mu\lra V_\nu$ such that
$\nu=\lambda+\mu$ and $\phi^\vee(x_\nu)=x_\lambda\otimes x_\mu$ we have
$\phi(\upsilon_\lambda\otimes \upsilon_\mu)=\upsilon_\nu$;

b) For any $\bG$-morphism $\phi:\ V_\lambda\otimes V_\mu\lra V_\nu$ such that
$\nu<\lambda+\mu$ we have $\phi(\upsilon_\lambda\otimes\upsilon_\mu)=0$;

c) For any $\bG$-morphism
$\varphi:\ V_\lambda^\vee\otimes V_\mu^\vee\lra V_\nu^\vee$ such that
$\nu=\lambda+\mu$ and $\varphi(x_\lambda\otimes x_\mu)=x_\nu$ we have
$\varphi(u_\lambda\otimes u_\mu)=u_\nu$;

d) For any $\bG$-morphism
$\varphi:\ V_\lambda^\vee\otimes V_\mu^\vee\lra V_\nu^\vee$ such that
$\nu<\lambda+\mu$ we have $\varphi(u_\lambda\otimes u_\mu)=0$;

e) $\langle u_\lambda,\upsilon_\lambda\rangle=1$.

{\em Proof.} We are considering the
systems $(v,u)=\left(v_\la\in V_\lambda,\ u_\la\in V_\lambda^\vee,
\ {\la\in X^+}\right)$ such that both $v$ and $u$ are Pl\"ucker sections
and $\langle v,u\rangle =1$, i.e., $\langle v_\la,u_\la\rangle =1$
for each $\la$.

These form a $\bG$-torsor and we have fixed its element $(y,x)$,
which we will use to think of this torsor as a Pl\"ucker picture
of $\bG$.

The stabilizers $\bG_v$ and $\bG_u$ are the
unipotent radicals of the opposite Borel subgroups, for instance
$\bG_x=\bN$ and $\bG_y=\bN_-$.
So this torsor canonicaly maps into the
open $\bG$-orbit in $\bX\tim\bX$ and the fiber at $(\bB',{\bB}'')$ is
a torsor for a Cartan subgroup $\bB'\cap\bB''$.
 $\Box$

\subsubsection{}
\label{suda}
Given a Pl\"ucker section $(v_\lambda)_{\lambda\in X^+}$, the collection
of meromorphic sections
$(x_\lambda\in\CV_\lambda^\vee,\frac{v_\lambda}{Q_\lambda}\in\CV_\lambda)$
evidently satisfies the conditions
a)--e) of the above Lemma, and hence defines a meromorphic function
$g:\ \BA^1\lra \bG$. Let us list the properties of this function.

a) By the definition ~\ref{drinf} of the isomorphism $\psi$,
since $g$ fixes the Pl\"ucker section $x$ the function
$g$ actually takes values in $\bN\subset \bG$;

b) The argument similar to that used in ~\ref{tuda} shows that $g$ can be
extended to $\BP^1$, is regular at $\infty$, and $g(\infty)=1\in \bN$
(since $\frac{v_\lambda}{Q_\lambda}(\infty)=y_\la$,
$g(\infty)$ stabilizes $x_\la$ and $y_\la$ so it lies
in $\bN\cap\bN_-$);

c) Let $D=\pi_\alpha(v_\lambda)$ (see ~\ref{pi}) be the $I$-colored divisor
supported at the roots of $Q_\lambda$. Then $g$ is regular on $\BP^1-D$.

\subsubsection{}
\label{sjuda}
We define $\xi(v_\lambda)=(D,\CT,\tau)\in\bZ^\alpha$ as follows:
$D=\pi_\alpha(v_\lambda);\ \CT$ is the trivial $\bG$-torsor; the section $\tau$
is given by the function $g$.
Let us describe the corresponding $\bbH$-torsor
$\CT_{\tau,\bB_-}$
with meromorphic
section
$\tau_{\bB_-}$.
To describe an $\bbH$-torsor $\fL$ with a section $s$
it suffices to describe the induced $\BC^*$-torsors $\fL_\lambda$
with sections $s_\lambda$ for all characters $\lambda:\ \bbH\lra\BC^*$.
In fact, it suffices to consider only
$\lambda\in X^+$.
Then $\fL_\lambda$ is given by the construction of
~\ref{obratno}, and $s_\lambda=\frac{v_\lambda}{Q_\lambda}$.

Thus, the conditions i),ii) of the Definition ~\ref{bZ} are evidently satisfied.

\subsubsection{}
\label{trivial}
To proceed with the inverse construction, we will need the following

{\em Lemma.} Suppose $(D,\CT,\tau)\in\bZ^\alpha$. Then $\CT$ is trivial
and has a canonical section $\varsigma$.

{\em Proof.} By the construction ~\ref{type}, $\CT$ is induced from the
$\bB$-torsor
$\barr{\bB\cdd \tau}$.
 By the Definition ~\ref{bZ}, the induced $\bbH$-torsor
$\CT_{\tau,\bB}$ is trivial, that is,
$\barr{\bB\cdd \tau}$
can be further reduced to an
$\bN$-torsor. But any $\bN$-torsor over $\BP^1$ is trivial since
$H^1(\BP^1,\bV)=0$ for any unipotent group $\bV$
(induction in the lower central series). $\Box$

\subsubsection{}
According to the above Lemma, we can find a unique section $\varsigma$ of
$\CT$ defined over the whole $\BP^1$ and such that $\varsigma(\infty)=
\tau(\infty)$. Hence a triple $(D,\CT,\tau)\in\bZ^\alpha$ canonically
defines a meromorphic function
$$
g\df\ \tau\varsigma^{-1}:\ \BP^1\lra \bG
,$$
i.e., $g(x)\cdd\varsigma(x)=\tau(x),\ x\in\BP^1$.
One sees immediately that $g$ enjoys the properties ~\ref{suda}a)--c).
Now
we can apply the Lemma ~\ref{drinf} in the opposite direction and obtain
from $g$ a collection
$\psi^{-1}(g)=(x_\lambda,\ti\upsilon_\lambda)_{\lambda\in X^+}$
with $\ti\upsilon_\la$ a certain meromorphic sections of $\CV_\la$.
According to ~\ref{polynom}, the divisor $D$ defines a collection
of unitary polynomials $(Q_\lambda)_{\lambda\in X^+}$, and we can define
$v_\lambda\df\ Q_\lambda\ti\upsilon_\lambda$.
One checks easily that $(v_\lambda)\in\sZ^\alpha$, and
$(D,\CT,\tau)=\xi(v_\lambda)$.

In particular, $\xi:\ \sZ^\alpha\lra\bZ^\alpha$ is an isomorphism.

\subsection{}
\label{Z}
We conclude that $\CZ^\alpha,\sZ^\alpha,\bZ^\alpha$ are all the same
and all  maps to $\BA^\al$ coincide.
We preserve the notation $\CZ^\alpha$ for this {\em Zastava}
space, and $\pi_\alpha$ for its
projection onto $\BA^\alpha$. We combine the properties ~\ref{dimension},
~\ref{affine},~\ref{factorization},~\ref{irred} into the following

{\bf Theorem.} a) $\CZ^\alpha$ is an irreducible affine algebraic variety
of dimension $2|\alpha|$;

b) For any decomposition $\alpha=\beta+\gamma,\ \beta,\gamma\in\BN[I]$, and
a pair of disjoint subsets $U,\Upsilon\subset\BA^1$, we have the
{\em factorization property}
(notations of ~\ref{bZU} and ~\ref{factorization}):
$$\CZ^{\beta,\gamma}_{U,\Upsilon}=\CZ^\beta_U\times\CZ^\gamma_\Upsilon$$

\bigskip

\centerline{\bf CHAPTER 2. The category $\PS$}

\section{Schubert stratification}

\subsection{}
\label{sig}
We will stratify $\ZZ^\al$ in stages.
We denote by $\CQ^\al\suppp\qp^\al\suppp \qc^\al $,
respectively the variety
of all quasimaps of degree $\al$ and the subvarieties of
the quasimaps defined at $0$ and of genuine maps.
In the same way we denote the varieties of based quasimaps
$\ZZ^\al\suppp\zp^\al\df\ZZ^\al\cap\qp^\al
\suppp
\zc^\al\df\ZZ^\al\cap \qc^\al=$ based maps of degree $\al$.

Recall (see ~\ref{sigma}) the map
$\sigma_{\beta,\ga}: \CQ^\be\tim C\gaa\ra \CQ^{\be+\ga},
\ \si_{\be,\ga}(f,D)=f(-D)$.
For $\be\le \al(=\be+\ga)$, it restricts
to the embedding
$ \CQ^\be\inj \CQ^{\al},
\ f\mm f(-(\al-\be)\cdd 0)=f(\ (\be-\al)\cdd 0)$, and in particular
$\CZ^\beta\hra\CZ^\alpha$.

\subsection{}
\label{coarse}
In the first step we stratify $\ZZ^\al$ according to the singularity at 
$0$.
It follows immediately
from the Theorem ~\ref{strat M} that
$$
\CZ^\alpha\cong
\bigsqcup_{0\leq\beta\leq\alpha}\dZ^\beta
.$$

The closed embedding of a stratum $\dZ^\beta$ into $\CZ^\alpha$ will be
denoted by $\sigma_{\beta,\alpha-\beta}$.

\subsection{}
Next, we stratify the quasimaps $\dZ^\al$ defined at $0$, according to the
singularity on $\ccs$.
Again, it follows immediately
from the Theorem ~\ref{strat M} that
$$
\zp^\alpha\cong
\bigsqcup_{0\leq\beta\leq\alpha}\zc^\beta\tim (\ccs)^{\al-\be}
.$$

\subsection{}
\label{Schubert}
One more refinement comes from the decomposition of the
flag variety $\bX$ into the $\bB$-invariant Schubert cells.
Given an element $w$ in the Weyl group $\CW_f$,
we define the locally closed subvarieties
({\em Schubert strata}) $\dZ^{\alpha}_w\subset\dZ^\alpha$
and $\zc^{\alpha}_w\subset\zc^\alpha$, as the sets of
quasimaps $f$ such that $f(0)\in\bX_w$.
The closure of $\dZ^{\alpha}_w$ in $\CZ^\alpha$ will be denoted by
$\oCZ^{\alpha}_w$.
Evidently, 
$$\dZ^\alpha=\bigsqcup_{w\in \CW_f}\dZ^{\alpha}_w
\aand \zc^\alpha=\bigsqcup_{w\in \CW_f}\zc^{\alpha}_w.$$
Beware that $\dZ^{\alpha}_w$ may happen to be empty: e.g. for $\alpha=0,
w\not=w_0$.

\subsubsection{}
\label{fineS Q}
Finally, the last refinement comes from the
diagonal stratification of the configuration space $(\ccs)^\delta=
\bigsqcup_{ \Ga\in\fP(\delta)}(\ccs)_\Ga^\delta$.
Altogether,
we obtain the following stratifications of $\ZZ^\al$:
$$
\CZ^\alpha\cong
\bigsqcup_{\alpha\ge \beta}\dZ^\beta
\ \ (\text{{\em coarse stratification}})
$$
$$
\cong
\bigsqcup^{\alpha\ge \beta\ge\ga}_{\Gamma\in\fP(\beta-\gamma)}
\zc^\ga\tim (\ccs)^{\be-\ga}_\Gamma \ \
(\text{{\em fine stratification}})
$$
$$
\cong\bigsqcup^{\alpha\ge \beta\ge\ga}_{ w\in \CW_f,\ \Ga\in\fP(\be-\ga)}
\zc^\ga_w \tim (\ccs)_\Ga^{\beta-\gamma}
\ \ (\text{{\em fine Schubert stratification}})
.$$
Similarly, we have the {\em fine stratification}
(resp. {\em fine Schubert stratification}) of $\CQ^\alpha$:
$$\CQ^\alpha=
\bigsqcup^{\alpha\ge \beta\ge\ga}_{\Gamma\in\fP(\beta-\gamma)}
\qc^\ga\tim (\BP^1-0)^{\be-\ga}_\Gamma=
\bigsqcup^{\alpha\ge \beta\ge\ga}_{ w\in \CW_f,\ \Ga\in\fP(\be-\ga)}
\qc^\ga_w \tim (\BP^1-0)_\Ga^{\beta-\gamma}$$
Here $\qc^\gamma_w\subset\qc^\gamma$ denotes the locally closed subspace of
maps $\BP^1\to\bX$ taking value in $\bX_w\subset\bX$ at $0\in\BP^1$.
The strata $\qc^\gamma_w\times(\BP^1-0)_\Gamma^{\beta-\gamma}$ are evidently
smooth.

Note that the strata $\zp^\al_w$
are not necessarily smooth in general, e.g. for
$\bG=SL_3,\ \alpha$ the sum of simple coroots, $w=w_0$.
To understand the ``fine Schubert strata''
$\zc^\ga_w \tim (\ccs)_\Ga^\beta$ we need to understand the varieties
$\zc^\ga_w$.

\subsection{Conjecture}
\label{conj}
For $\gamma\in\BN[I],w\in\CW_f$ the variety
$\zc^\gamma_w$ is smooth. Hence the ``fine Schubert stratification''
is really a stratification.

\subsubsection{Lemma}
\label{cheap}
For $\gamma$ sufficiently dominant (i.e. $\langle\gamma,i'\rangle>10$) and
arbitrary $w\in\CW_f$ the variety $\zc^\gamma_w$ is smooth.

{\em Proof.} Let us consider the map
$\varrho_\gamma:\ \qc^\gamma\lra\bX\times\bX,\ f\mapsto(f(0),f(\infty))$.
We have $\oZ^\gamma=\varrho_\gamma^{-1}(\bX_w,\bB_-)$.
It suffices to prove that $\varrho_\gamma$ is smooth and surjective.
Recall that the tangent
space $\Theta_f$ to $\qc^\gamma$ at $f\in\qc^\gamma$ is canonically
isomorphic to $H^0(\BP^1,f^*\CT\bX)$. Let us interpret $\bX$ as a variety of
Borel subalgebras of $\fg$. We denote $f(0)$ by $\fb_0$, and $f(\infty)$ by
$\fb_\infty$. So we have to prove that the canonical map
$H^0(\BP^1,f^*\CT\bX)\lra\CT_{\fb_0}\bX\oplus\CT_{\fb_\infty}\bX$ is surjective.
To this end it is enough to have $H^1(\BP^1,f^*\CT\bX(-0-\infty))=0$.
This in turn holds whenever $\gamma$ is sufficiently dominant.
$\Box$

\subsubsection{Lemma}
\label{dim Schubert}
For $\gamma$ sufficiently dominant we have
$\dim\oZ^\gamma_w=2|\gamma|-\dim\bX+\dim\bX_w$.

{\em Proof.} The same as the proof of ~\ref{cheap}. $\Box$

\subsubsection{Remark} Unfortunately, one cannot prove the conjecture
~\ref{conj} for arbitrary $\gamma$ the same way
as the Lemma ~\ref{cheap}: for arbitrary
$\gamma$ the map $\varrho_\gamma$ is not smooth. The simplest example
occurs for $\bG=SL_4$ when $\gamma$ is twice the sum of simple coroots.
This example was found by A.Kuznetsov.

\section{Factorization}

This section follows closely \S4 of ~\cite{fs}.

\subsection{}
\label{fake}
Now we replace the maps into the flag variety $\bX$ with the maps
from $\BP^1$ to the product $\bX\tim Y=\sqcup_{\chi\in Y}\ \bX_\chi$.
So for arbitrary  $\chi\in Y$ and $\alpha\in\BN[I]$
we obtain the spaces $\CZ^\alpha_\chi$ of based maps into $\bX_\chi$
and it makes sense now  
to add the  subscript $\chi$ to all the
strata (coarse, Schubert, fine) defined in the previous section.

We will consider a system $\ZZ$ of  varieties
$\CZ^\alpha_\chi,\ \al,\ga\in Y$, together with  two kinds of maps
defined for any $\beta,\ga\in\BN[I]$:

a) closed embeddings,
$$\sigma^{\beta,\gamma}_\chi:\
\CZ^\beta_\chi\hra\CZ^{\beta+\gamma}_{\chi+\gamma}
,$$

b) factorization identifications
$$\CZ^{\beta,\gamma}_{\chi,\Ue,\Upe}=
\CZ^\beta_{\chi,\Ue}\times\CZ^\gamma_{\chi-\beta,\Upe}$$
defined for $\varepsilon>0$ and
$U_\varepsilon\df\ 
\{z\in\BC,\ |z|<\varepsilon\}$, and
$\Upsilon_\varepsilon\df\
\{z\in\BC,\ |z|>\varepsilon\}$.

Of course, without the subscript these are the
factorizations from  ~\ref{factorization}
and the  embeddings
from  ~\ref{sig}.

\subsection{}
\label{snop}
We will denote by $\IC^\alpha_\chi$ the perverse $IC$-extension of the
constant sheaf at the generic point of $\CZ^\alpha_\chi$.

The following definition makes sense only modulo the validity of conjecture
~\ref{conj}.

{\bf Definition.} A {\em snop} $\CK$ is the following collection of data:

a) $\chi=\chi(\CK)\in Y$, called the {\em support estimate} of $\CK$;

b) For any $\alpha\in\BN[I]$, a perverse sheaf $\CK^\alpha_\chi$ on
$\CZ^\alpha_\chi$ smooth along the fine Schubert stratification;

c) For any $\beta,\gamma\in\BN[I],\ \varepsilon>0$,
a {\em factorization isomorphism}
$$
\CK^{\beta+\gamma}_\chi|_{\CZ^{\beta,\gamma}_{\chi,\Ue,\Upe}}
\iso
\CK^\beta_\chi|_{\CZ^\beta_{\chi,\Ue}}
\boxtimes
\IC^\gamma_{\chi-\beta}|_{\CZ^\gamma_{\chi-\beta,\Upe}}$$
satisfying the {\em associativity constraints} as in ~\cite{fs}, \S\S 3,4.
We spare the reader the explicit formulation of these constraints.

\subsection{}
\label{awful}
Since at the moment the conjecture ~\ref{conj} is unavailable we will provide
an ugly provisional substitute of the Definition ~\ref{snop}. Namely, recall
that $\CZ^\alpha=
\sqcup_{\alpha\geq\beta\geq\gamma}\oZ^\gamma\times(\BC^*)^{\beta-\gamma}$.
We introduce an open subvariety
$$\ddZ^\alpha=\bigsqcup_{\alpha\geq\beta\geq\gamma\gg0}\oZ^\gamma\times
(\BC^*)^{\beta-\gamma}$$
The union is taken over sufficiently dominant $\gamma$, i.e. such that
$\langle\gamma,i'\rangle>10$ for any $i\in I$. Certainly, if $\alpha$ itself
is not sufficiently dominant, $\ddZ^\alpha$ may happen to be empty.
We have the fine Schubert stratification
$$\ddZ^\alpha=\bigsqcup^{\alpha\geq\beta\geq
\gamma\gg0}_{w\in\CW_f,\ \Gamma\in\fP(\beta-\gamma)}\oZ^\gamma_w\times
(\BC^*)^{\beta-\gamma}_\Gamma$$ with smooth strata (see the Lemma ~\ref{cheap}).

Now we can repeat the Definition ~\ref{snop} replacing $\CZ^\alpha_\chi$
by $\ddZ^\alpha_\chi$. Thus in ~\ref{snop} b) we have to restrict ourselves
to sufficiently dominant $\alpha$, and in ~\ref{snop} c) $\beta$ has to be
sufficiently dominant as well.

\subsubsection{}
In what follows we use the Definition ~\ref{snop}. The reader unwilling to
believe in the Conjecture ~\ref{conj} will readily substitute the conjectural
Definition ~\ref{snop} with the provisional working Definition ~\ref{awful}.

\subsection{Examples}
We define the {\em irreducible} and {\em standard} snops.

\subsubsection{}
\label{CL}
Let us describe a snop $\CL(w,\chi)$ for $\chi\in Y,\ w\in \CW_f$.

a) The support of $\CL(w,\chi)$ is $\chi$.

b) $\CL(w,\chi)^\alpha_\chi$ is the irreducible $IC$-extension
$\IC(\oCZ^{\alpha}_{w,\chi})=j_{!*}\IC(\dZ^{\alpha}_{w,\chi})$
of the perverse $IC$-sheaf on the Schubert stratum $\dZ^{\alpha}_{w,\chi}\subset
\dZ^\alpha_\chi$. Here $j$ stands for the affine open embedding
$\dZ^{\alpha}_{w,\chi}\hra\oCZ^{\alpha}_w$.

In particular, $\IC(\oCZ^\alpha_{w_0,\chi})=\IC^\alpha_\chi$.


c) Evidently, $\oCZ^{\beta}_{w,\chi,\Ue}$
(resp. $\oCZ^{\beta,\gamma}_{w,\chi,\Ue,\Upe}$) is open in
$\oCZ^{\beta}_{w,\chi}$ (resp. $\oCZ^{\alpha}_{w,\chi}$) for any
$\beta+\gamma=\alpha$. Moreover, $\oCZ^{\beta,\gamma}_{w,\chi,\Ue,\Upe}
=\oCZ^{\beta}_{w,\chi,\Ue}\times\CZ^\gamma_{\chi-\beta,\Upe}$.
This induces the desired factorization isomorphism.

\subsubsection{}
\label{CM}
If we replace in ~\ref{CL}b) above $j_{!*}\IC(\dZ^{\alpha}_{w,\chi})$ by
$j_!\IC(\dZ^{\alpha}_{w,\chi})=:\CM(w,\chi)^\alpha_\chi$
(resp. $j_*\IC(\dZ^{\alpha}_{w,\chi})=:\CalD\CM(w,\chi)^\alpha_\chi$)
we obtain the snop $\CM(w,\chi)$ (resp. $\CalD\CM(w,\chi)$).

\subsection{}
Given a snop $\CK$ with support $\chi$, and $\eta\geq\chi,\
\alpha\in\BN[I]$, we define a sheaf $'\CK^\alpha_\eta$ on $\CZ^\alpha_\chi$
as follows. We set $\gamma\df\ \eta-\chi$. If $\alpha\geq\gamma$ we set
$$'\CK^\alpha_\eta\df\ (\sigma^{\alpha-\gamma,\gamma}_\chi)_*
\CK^{\alpha-\gamma}_\chi$$
(for the definition of $\sigma$ see ~\ref{fake}).
Otherwise we set $'\CK^\alpha_\eta\df\ 0$.

It is easy to see that the factorization isomorphisms for $\CK$ induce
similar isomorphisms for $'\CK$, and thus we obtain a snop $'\CK$ with support
$\eta\geq\chi$.

\subsection{}
We define the category $\widetilde\PS$ of snops.

\subsubsection{}
\label{morphism}
Given two snops $\CF,\CK$ we will define the morphisms $\Hom(\CF,\CK)$
as follows. Let $\eta\in Y$ be such that $\eta\geq\chi(\CF),\chi(\CK)$.
For $\alpha=\beta+\gamma\in\BN[I]$ we consider the following composition:
$$\vartheta^{\beta,\gamma}_\eta:\ \Hom_{\CZ^\alpha_\eta}('\CF^\alpha_\eta,
'\CK^\alpha_\eta)\lra\Hom_{\CZ^{\beta,\gamma}_{\Ue,\Upe}}
('\CF^\alpha_\eta|_{\CZ^{\beta,\gamma}_{\Ue,\Upe}},
'\CK^\alpha_\eta|_{\CZ^{\beta,\gamma}_{\Ue,\Upe}})\iso$$
$$\Hom_{\CZ^\beta_{\eta,\Ue}\times\CZ^\gamma_{\eta-\beta,\Upe}}
('\CF^\beta_\eta|_{\CZ^\beta_{\eta,\Ue}}\boxtimes\IC^\gamma_{\eta-\beta}
|_{\CZ^\gamma_{\eta-\beta,\Upe}},
'\CK^\beta_\eta|_{\CZ^\beta_{\eta,\Ue}}\boxtimes\IC^\gamma_{\eta-\beta}
|_{\CZ^\gamma_{\eta-\beta,\Upe}})=
\Hom_{\CZ^\beta_{\eta,\Ue}}('\CF^\beta_\eta|_{\CZ^\beta_{\eta,\Ue}},
'\CK^\beta_\eta|_{\CZ^\beta_{\eta,\Ue}})$$
(the second isomorphism is induced by the factorization isomorphisms for
$'\CF$ and $'\CK$, and the third equality is just K\"unneth formula).

Now we define
$$\Hom(\CF,\CK)\df\ \dirlim_\eta\invlim_\alpha
\Hom_{\CZ^\alpha_\eta}('\CF^\alpha_\eta,'\CK^\alpha_\eta)$$

Here the inverse limit is taken over $\alpha\in\BN[I]$, the transition
maps being $\vartheta^{\beta,\alpha-\beta}_\eta$, and the direct limit is
taken over $\eta\in Y$ such that $\eta\geq\chi(\CF),\chi(\CK)$, the
transition maps being induced by the obvious isomorphisms
$\Hom_{\CZ^\alpha_\eta}('\CF^\alpha_\eta,'\CK^\alpha_\eta)=
\Hom_{\CZ^{\alpha+\gamma}_{\eta+\gamma}}('\CF^{\alpha+\gamma}_{\eta+\gamma},
'\CK^{\alpha+\gamma}_{\eta+\gamma})$.

\subsubsection{}
With the above definition of morphisms and obvious composition, the snops
form a category which we will denote by $\widetilde\PS$.

\subsection{}
\label{PS}
Evidently, the snops $\CL(w,\chi)$ are irreducible objects of $\widetilde\PS$.
It is easy to see that any irreducible object of $\widetilde\PS$ is isomorphic
to some $\CL(w,\chi)$.

We define the category $\PS$ of {\em finite snops}
as the full subcategory of $\widetilde\PS$ formed
by the snops of finite length. It is an abelian category. We will see later
that $\CM(w,\chi)$ and $\CalD\CM(w,\chi)$ (see ~\ref{CM}) lie in $\PS$ for any
$w,\chi$.

One can prove the following very useful technical lemma exactly as in
~\cite{fs}, 4.7.

\subsubsection{Lemma}
\label{stabilize}
Let $\CF,\CK$ be two finite snops. Let $\eta\geq\chi(\CF),\chi(\CK)$.
There exists $\beta\in\BN[I]$ such that for any $\alpha\geq\beta$
the canonical maps $\Hom(\CF,\CK)\lra
\Hom_{\CZ^\alpha_\eta}('\CF^\alpha_\eta,'\CK^\alpha_\eta)$
are all isomorphisms. $\Box$


\bigskip

\centerline{\bf CHAPTER 3. Convolution with affine Grassmannian}

\section{Pl\"ucker model of affine Grassmannian}

\subsection{}
Let $\CG$ be the usual affine Grassmannian $\bG((z))/\bG[[z]]$. It is the
ind-scheme representing the functor of isomorphism classes of pairs
$(\CT,\tau)$ where $\CT$ is a $\bG$-torsor on $\BP^1$, and $\tau$ is its
section (trivialization) defined off $0$ (see e.g. ~\cite{mv}). It is equipped
with a natural action of proalgebraic group $\bG[[z]]$, and we are going
to describe the orbits of this action. It is known (see e.g. {\em loc. cit.})
that these orbits are numbered by dominant cocharacters $\eta\in Y^+\subset Y$.

Here $Y^+\subset Y$ stands for the set of cocharacters $\eta$ such that
$\langle\eta,i'\rangle\geq0$ for any $i\in I$. For $\eta\in Y^+$ we denote
the corresponding $\bG[[z]]$-orbit in $\CG$ by $\CG_\eta$, and we denote
its closure by $\oCG_\eta$.

Recall that for a dominant character $\lambda\in X^+$ we denote by $V_\lambda$
the corresponding irreducible $\bG$-module, and we denote by $\CV_\lambda$
the trivial vector bundle $V_\lambda\otimes\CO_{\BP^1}$ on $\BP^1$.

\subsection{Proposition}
\label{closure}
The orbit closure $\oCG_\eta\subset\CG$ is the space of collections
$(\CU_\lambda)_{\lambda\in X^+}$ of vector bundles on $\BP^1$ such that

a) $\CV_\lambda(-\langle\eta,\lambda\rangle0)\subset\CU_\lambda\subset
\CV_\lambda(\langle\eta,\lambda\rangle0)$, or equivalently,
$\CU_\lambda(-\langle\eta,\lambda\rangle0)\subset\CV_\lambda\subset
\CU_\lambda(\langle\eta,\lambda\rangle0)$;

b) $\deg\CU_\lambda=\deg\CV_\lambda=0$, or in other words,
$\dim\CV_\lambda(\langle\eta,\lambda\rangle0)/\CU_\lambda=
\langle\eta,\lambda\rangle\dim V_\lambda$;

c) For any surjective $G$-morphism $\phi:\ V_\lambda\otimes V_\mu\lra V_\nu$
and the corresponding morphism $\phi:\ \CV_\lambda\otimes\CV_\mu\lra\CV_\nu$
(hence $\phi:\ \CV_\lambda(\langle\eta,\lambda\rangle0)\otimes
\CV_\mu(\langle\eta,\mu\rangle0)\lra
\CV_\nu(\langle\eta,\lambda+\mu\rangle0)$) we have
$\phi(\CU_\lambda\otimes\CU_\mu)=\CU_\nu$.

{\em Proof.} $\bG$-torsor on a curve $C$ is the same as a tensor functor from
the category of $\bG$-modules to the category of vector bundles on $C$. $\Box$

\subsection{}
\label{loc}
Let us give a local version of the above Proposition. Recall that
$\CO=\BC[[z]]\subset\CK=\BC((z))$. For a finite-dimensional vector space $V$,
a {\em lattice} $\fV$ in $V\otimes\CK$ is an $\CO$-submodule of $V\otimes\CK$
{\em commeasurable} with $V\otimes\CO$, that is, such that
$(V\otimes\CO)\cap \fV$
is of finite codimension in both $V\otimes\CO$ and $\fV$.

{\bf Proposition.}
The orbit closure $\oCG_\eta\subset\CG$ is the space of collections
$(\fV_\lambda)_{\lambda\in X^+}$ of lattices in $V_\lambda\otimes\CK$ such that

a) $z^{\langle\eta,\lambda\rangle}(V_\lambda\otimes\CO)\subset
\fV_\lambda\subset
z^{-\langle\eta,\lambda\rangle}(V_\lambda\otimes\CO)$, or equivalently,
$z^{\langle\eta,\lambda\rangle}\fV_\lambda\subset V_\lambda\otimes\CO\subset
z^{-\langle\eta,\lambda\rangle}\fV_\lambda$;

b) $\dim(z^{-\langle\eta,\lambda\rangle}(V_\lambda\otimes\CO)/\fV_\lambda)=
\langle\eta,\lambda\rangle\dim V_\lambda$;

c) For any surjective $G$-morphism $\phi:\ V_\lambda\otimes V_\mu\lra V_\nu$
and the corresponding morphism
$\phi:\ (V_\lambda\otimes\CO)\otimes(V_\mu\otimes\CO)\lra(V_\nu\otimes\CO)$
(hence $\phi:\ z^{-\langle\eta,\lambda\rangle}(V_\lambda\otimes\CO)\otimes
z^{-\langle\eta,\mu\rangle}(V_\mu\otimes\CO)\lra
z^{-\langle\eta,\lambda+\mu\rangle}(V_\nu\otimes\CO)$), we have
$\phi(\fV_\lambda\otimes \fV_\mu)=\fV_\nu$.

$\Box$

\subsection{}
\label{Iwahori}
Let $\bI\subset\bG[[z]]$ be the Iwahori subgroup; it is formed
by all $g(z)\in\bG[[z]]$ such that $g(0)\in\bB\subset\bG$. We will denote by
$\CP(\CG,\bI)$ the category of perverse sheaves on $\CG$
with finite-dimensional support, constant along $\bI$-orbits.
The stratification of $\CG$ by $\bI$-orbits is a certain refinement of
the stratification $\CG=\sqcup_{\eta\in Y^+}\CG_\eta$. Namely, each $\CG_\eta$
decomposes into $\bI$-orbits numbered by $\CW_f/\CW_\eta$ where $\CW_\eta$
stands for the stabilizer of $\eta$ in $\CW_f$. For $w\in\CW_f/\CW_\eta$
we will denote the corresponding $\bI$-orbit by $\CG_{w,\eta}$. Let us
introduce a Pl\"ucker model of $\CG_{w,\eta}$.

\subsubsection{} For $\eta\in Y^+$ let $I_\eta\subset I$ be the set of all $i$
such that $\langle\eta,i'\rangle=0$ (thus for $i\not\in I_\eta$ we have
$\langle\eta,i'\rangle>0$). Then $\CW_\eta$ is generated by the simple
reflections $\{s_i, i\in I_\eta\}$. Let $\bP(I_\eta)$ be the corresponding
parabolic subgroup (e.g. for $I_\eta=\emptyset$ we have $\bP(I_\eta)=\bB$,
while for $I_\eta=I$ we have $\bP(I_\eta)=\bG$). Let $\bX(I_\eta)=
\bG/\bP(I_\eta)$ be the corresponding partial flag variety. The $\bB$-orbits
on $\bX(I_\eta)$ are naturally numbered by $\CW_f/\CW_\eta:\ \bX(I_\eta)=
\sqcup_{w\in\CW_f/\CW_\eta}\bX(I_\eta)_w$. The Pl\"ucker embedding realizes
$\bX$ as a closed subvariety in $\prod_{i\in I}\BP(V_{\omega_i})$. Its image
under the projection
$\prod_{i\in I}\BP(V_{\omega_i})\lra\prod_{i\not\in I_\eta}
\BP(V_{\omega_i})$ exactly coincides with $\bX(I_\eta)$.

\subsubsection{Lemma-Definition}



a) For $\eta=\sum_{i\in I}n_ii$,
and $(\CU_\lambda)_{\lambda\in X^+}\in\CG_\eta$ we have\\
$\dim(\CU_{\omega_i}+
\CV_{\omega_i}((n_i-1)\cdot0)/\CV_{\omega_i}((n_i-1)\cdot0))=
\dim V_{\omega_i}^{{\bf U}(I_\eta)}$ where ${\bf U}(I_\eta)$ is the unipotent
radical of $\bP(I_\eta)$.

b) Thus $\CU_{\omega_i}, i\in I,$ defines a subspace $K_i$ in
$\CV_{\omega_i}(n_i\cdot0)/\CV_{\omega_i}((n_i-1)\cdot0)=
V_{\omega_i}$. This collection
of subspaces $(K_i)_{i\in I}\in\prod_{i\in I}\operatorname{Gr}(V_{\omega_i})$
satisfies the Pl\"ucker relations and thus gives a point in $\bX(I_\eta)$;

c) We will denote by $\br$ the map $\CG_\eta\lra\bX(I_\eta)$ defined in b);

d) For $w\in\CW_f/\CW_\eta$ we have $\CG_{w,\eta}=\br^{-1}(\bX(I_\eta)_w)$.
$\Box$

We are obliged to D.Gaitsgory who pointed out a mistake in the earlier
version of the above Lemma.

\subsubsection{}
For $\theta\in Y$ we consider the corresponding homomorphism $\theta:\
\BC^*\lra\bH\subset\bG$ as a formal loop $\theta(z)\in\bG((z))$.
It projects to the same named point $\theta(z)\in\bG((z))/\bG[[z]]=\CG$.
There is a natural bijection between the set of $\theta(z),\ \theta\in Y$,
and the set of Iwahori orbits: each Iwahori orbit $\CG_{w,\eta}$ contains
exactly one of the above points, namely, the point $w\eta(z)$.

\subsection{}
\label{PBD}
Recall the Beilinson-Drinfeld avatar $\bZ^\alpha$ of the Zastava space
$\CZ^\alpha$ (see ~\ref{bZ}). In this subsection we will give a Pl\"ucker model
of $\bZ^\alpha$.

{\bf Proposition.} $\bZ^\alpha$ is the the space of pairs
$(D,(\fU_\lambda)_{\lambda\in X^+})$ where $D\in\BA^\alpha$ is an
$I$-colored effective divisor, and
$(\fU_\lambda)_{\lambda\in X^+}$ is a collection of vector bundles on $\BP^1$
such that

a) $\CV_\lambda(-\infty D)\subset\fU_\lambda\subset\CV_\lambda(+\infty D)$;

b) $\CV_\lambda^\bN\subset\fU_\lambda\supset
\CV_\lambda^{\bN_-}(-\langle D,\lambda\rangle)$ (notations of ~\ref{sigma}),
the first inclusion being a {\em line subbundle} (and the second an invertible
subsheaf);

c) $\deg\fU_\lambda=0$;

d) For any surjective $G$-morphism $\phi:\ V_\lambda\otimes V_\mu\lra V_\nu$
and the corresponding morphism $\phi:\ \CV_\lambda\otimes\CV_\mu\lra\CV_\nu$
(hence $\phi:\ \CV_\lambda(+\infty D)\otimes
\CV_\mu(+\infty D)\lra
\CV_\nu(+\infty D)$) we have
$\phi(\fU_\lambda\otimes\fU_\mu)=\fU_\nu$.

{\em Proof.} Obvious. $\Box$

\subsubsection{Remark}
\label{?}
Recall the isomorphism $\varpi^{-1}\xi:\ \bZ^\alpha\iso\CZ^\alpha$ constructed
in section 7. Let us describe it in terms of ~\ref{PBD}. The Lemma
~\ref{trivial} says that there is a unique system of isomorphisms
$\iota_\lambda:\ \fU_\lambda\iso\CV_\lambda,\ \lambda\in X^+$, identical at
$\infty\in\BP^1$ and compatible with tensor multiplication.
Then $\varpi^{-1}\xi(D,(\fU_\lambda)_{\lambda\in X^+})=(\fL_\lambda\subset
\CV_\lambda)_{\lambda\in X^+}$ where $\fL_\lambda=
\iota_\lambda(\CV_\lambda^{\bN_-}(-\langle D,\lambda\rangle))$.

\subsection{}
\label{stack}
Let $\fM$ be the scheme representing the
functor of isomorphism classes of $\bG$-torsors on $\BP^1$ equipped with
trivialization in the formal neighbourhood of $\infty\in\BP^1$
(see ~\cite{ka} and ~\cite{kt1}).

\subsubsection{}
\label{strat stack}
The scheme $\fM$ is stratified by the locally closed subschemes $\fM_\eta:\
\fM=\sqcup_{\eta\in Y^+}\fM_\eta$ according to the isomorphism types of
$\bG$-torsors. Namely, due to Riemann's classification, for a $\bG$-torsor
$\CT$ and any
$\lambda\in X^+$ the associated vector bundle $\CV_\lambda^\CT$ decomposes
as a direct sum of line bundles $\CO(r_k^\lambda)$ of well-defined degrees
$r_1^\lambda\geq\ldots\geq r^\lambda_{\dim V_\lambda}$. Then $\CT$ lies in the
stratum $\fM_\eta$ iff $r_1^\lambda=\langle\eta,\lambda\rangle$.

For any $\eta\in Y^+$ the union of strata $\fM^\eta:=
\sqcup_{Y^+\ni\chi\leq\eta}\fM_\chi$ forms an open subscheme of $\fM$.
This subscheme is a projective limit of schemes of finite type, all the maps
in projective system being fibrations with affine fibers. Moreover,
$\fM^\eta$ is equipped with a free action of a prounipotent group $\bG^\eta$
(a congruence subgroup in $\bG[[z^{-1}]]$) such that the quotient $\ufM^\eta$
is a smooth scheme of finite type. The theory of perverse sheaves on $\fM$
smooth along the stratification by $\fM_\eta$ is developed in ~\cite{kt1}.
We will refer the reader to this work, and will freely use such perverse
sheaves, e.g. $\IC(\fM_\eta)$.

\subsubsection{}
\label{thin}
Restricting a trivialization of a $\bG$-torsor from $\BP^1-0$ to
the formal neighbourhood of $\infty\in\BP^1$ we obtain the closed
embedding $\bi:\ \CG\hra\fM$. The intersection of $\fM_\eta$ and $\CG_\chi$
is nonempty iff $\eta\leq\chi$, and then it is transversal. Thus, $\oCG_\eta
\subset\fM^\eta$. According to ~\cite{kt}, the composition $\oCG_\eta
\hookrightarrow\fM^\eta\lra\ufM^\eta$ is a closed embedding.

\subsubsection{}
\label{mapstack}
For a $\bG$-torsor $\CT$ and an irreducible $\bG$-module $V_\lambda$ we
denote by $\CV_\lambda^\CT$ the associated vector bundle. Following
~\ref{maps} and ~\ref{quasimaps} we define
for {\em arbitrary} $\alpha\in Y$ the scheme $\ofQ^\alpha$
(resp. $\fQ^\alpha$) representing the functor of isomorphism classes of pairs
$(\CT,(\fL_\lambda)_{\lambda\in X^+})$
where $\CT$ is a $\bG$-torsor trivialized in the formal neighbourhood of
$\infty\in\BP^1$, and $\fL_\lambda\subset\CV_\lambda^\CT,\
\lambda\in X^+$, is a collection of line subbundles (resp. invertible
subsheaves) of degree $\langle-\alpha,\lambda\rangle$ satisfying the Pl\"ucker
conditions (cf. {\em loc. cit.}). The evident projection $\ofQ^\alpha\lra\fM$
(resp. $\fQ^\alpha\lra\fM$) will be denoted by $\obp$ (resp. $\bp$).
The open embedding $\ofQ^\alpha\hra\fQ^\alpha$ will be denoted by $\bj$.
Clearly, $\bp$ is projective, and $\obp=\bp\circ\bj$.


The free action of prounipotent group $\bG^\eta$ on $\fM^\eta$ lifts to
the free action of $\bG^\eta$ on the open subscheme
$\bp^{-1}(\fM^\eta)\subset\fQ^\alpha$.
The quotient is a scheme of finite type $\ufQ^{\alpha,\eta}$ equipped with
the projective morphism $\bp$ to $\ufM^\eta$. There exists a
$\bG[[z^{-1}]]$-invariant stratification $\fS$ of $\fQ^\alpha$ such that
$\bp$ is stratified with respect to $\fS$ and the stratification
$\fM=\sqcup_{\eta\in Y^+}\fM_\eta$. One can define perverse sheaves on
$\fQ^\alpha$ smooth along $\fS$ following the lines of ~\cite{kt1}.
In particular, we have the irreducible Goresky-Macpherson sheaf
$\IC(\fQ^\alpha)$.

Following ~\ref{strat M} we introduce a decomposition of $\fQ^\alpha$
into a disjoint union of
locally closed subschemes according to the isomorphism types of $\bG$-torsors
and defects of invertible subsheaves:
$$\fQ^\alpha=\bigsqcup^{\eta\in Y^+}_{\beta\leq\alpha}\ofQ^\beta_\eta
\times C^{\alpha-\beta}$$ where $C=\BP^1$ and $\ofQ^\beta_\eta=
\obp^{-1}(\fM_\eta)\subset\ofQ^\beta$.

\subsection{IC sheaves}

\subsubsection{}
\label{GMZ}
The Goresky-MacPherson sheaf $\IC^\alpha$ on $\CZ^\alpha$
is smooth along stratification
$$\CZ^\alpha=
\bigsqcup^{\alpha\geq\beta\geq\gamma\geq0}_{\Gamma\in\fP(\beta-\gamma)}
\zc^\gamma\times(\BC^*)^{\beta-\gamma}_\Gamma$$
(cf. ~\ref{Schubert}). It is evidently constant along strata, so its stalk
at a point in $\zc^\gamma\times(\BC^*)^{\beta-\gamma}_\Gamma$
depends on the stratum only. Moreover, due to factorization property, it
depends not on $\alpha\geq\beta$ but only on their difference $\alpha-\beta\in
\BN[I]$. We will denote it by $\IC^{\alpha-\beta}_\Gamma$.
In case $\bG=SL_n$ these stalks were computed in ~\cite{ku}.

\subsubsection{}
\label{GMQ}
Recall (see ~\ref{strat M}) that $\CQ^\beta,\ \beta\in\BN[I]$, is stratified
by the type of defect:
$$\CQ^\beta=\bigsqcup^{\beta\geq\gamma\geq0}_{\Gamma\in\fP(\beta-\gamma)}
\oQ^\gamma\times C^{\beta-\gamma}_\Gamma$$
The Goresky-Macpherson sheaf $\IC(\CQ^\beta)$ on $\CQ^\beta$ is constant
along the strata. It is immediate to see that its stalk at any point in the
stratum $\oQ^\gamma\times C^{\beta-\gamma}_\Gamma$ is isomorphic, up to a
shift, to $\IC^0_\Gamma$. In particular, it depends on the defect only.

\subsubsection{}
\label{GM}
{\bf Proposition.}
a) The Goresky-Macpherson sheaf $\IC(\fQ^\beta)$ on $\fQ^\beta,\
\beta\in Y$, is constant along the locally closed subschemes
$$\fQ^\beta=\bigsqcup^{\beta\geq\gamma}_{\Gamma\in\fP(\beta-\gamma)}
\ofQ^\gamma\times C^{\beta-\gamma}_\Gamma$$

b) The stalk of $\IC(\fQ^\beta)$ at any point in the
$\ofQ^\gamma\times C^{\beta-\gamma}_\Gamma$ is isomorphic,
up to a shift, to $\IC^0_\Gamma$.

{\em Proof.} Will be given in ~\ref{later}. $\Box$

\subsubsection{}
\label{parity}
Let $\phi\in\fQ^\beta$. The stalk $\IC(\fQ^\beta)_\phi$ is a graded vector
space.

{\bf Conjecture.} {\em (Parity vanishing)} Nonzero graded parts of
$\IC(\fQ^\beta)_\phi$ appear in cohomological degrees of the same parity.

\subsubsection{Remark} In case $\bG=SL_n$ the conjecture follows from
the Proposition ~\ref{GM} and ~\cite{ku} 2.5.2. In general case the conjecture
follows from the unpublished results of G.Lusztig. We plan to prove it in
the next part.

\section{Convolution diagram}

\subsection{Definition} For $\alpha\in Y,\eta\in Y^+$ we define the
{\em convolution diagram} $\CG\CQ_\eta^\alpha$ as the space of collections
$(\CU_\lambda,\fL_\lambda)_{\lambda\in X^+}$ of vector bundles with invertible
subsheaves such that

a) $(\CU_\lambda)_{\lambda\in X^+}\in\oCG_\eta$, or in other words,
$(\CU_\lambda)_{\lambda\in X^+}$ satisfies the conditions ~\ref{closure} a)-c);

b) $\fL_\lambda\subset\CU_\lambda$ has degree $-\langle\alpha,\lambda\rangle$;

c) For any surjective $\bG$-morphism $\phi:\ V_\lambda\otimes V_\mu\lra V_\nu$
such that $\nu=\lambda+\mu$ we have (cf. ~\ref{closure} c)
$\phi(\CL_\lambda\otimes\CL_\mu)=\CL_\nu$;

d) For any $\bG$-morphism $\phi:\ V_\lambda\otimes V_\mu\lra V_\nu$ such that
$\nu<\lambda+\mu$ we have $\phi(\CL_\lambda\otimes\CL_\mu)=0$.

\subsubsection{} Let us denote by $\oGQ_\eta^\alpha$ the open subvariety in
$\CG\CQ_\eta^\alpha$ formed by all the collections $(\CU_\lambda,\fL_\lambda)$
such that $\fL_\lambda$ is a {\em line subbundle} in $\CU_\lambda$ for
any $\lambda\in X^+$. The open embedding
$\oGQ_\eta^\alpha\hra\CG\CQ_\eta^\alpha$ will be denoted by ~$\bj$.

\subsection{Definition} a) We define the projection $\bp:\ \CG\CQ_\eta^\alpha
\lra\oCG_\eta$ as $\bp(\CU_\lambda,\fL_\lambda)=(\CU_\lambda)$;

b) We define the map $\bq:\ \CG\CQ_\eta^\alpha\lra\CQ^{\alpha+\eta}$ as follows:
$$\bq(\CU_\lambda,\fL_\lambda)=(\fL_\lambda(-\langle\eta,\lambda\rangle0)
\subset\CU_\lambda(-\langle\eta,\lambda\rangle0)\subset\CV_\lambda)$$
(cf. ~\ref{closure} a) and the Definition ~\ref{quasimaps}).

\subsubsection{}
\label{some}
We will denote by $\obp$ the restriction of $\bp$ to the open subvariety
$\oGQ_\eta^\alpha\stackrel{\bj}{\hookrightarrow}\CG\CQ_\eta^\alpha$.

\subsubsection{Remark}
Note that $\CG\CQ_0^\alpha=\CQ^\alpha,\ \oGQ_0^\alpha=\qc^\alpha,\ \bq=\id$,
and $\bp$ is the projection to the point $\oCG_0$.

Note also that while in the Definition ~\ref{quasimaps} we imposed the
positivity condition $\alpha\in\BN[I]$ (otherwise $\CQ^\alpha$ would be empty)
here we allow arbitrary $\alpha\in Y$. It is easy to see that
$\CG\CQ_\eta^\alpha$ (as well as $\oGQ_\eta^\alpha$)
is nonempty iff $\eta+\alpha\in\BN[I]$.

\subsection{}
Let us give a local version of the convolution diagram. Recall the notations
of ~\ref{loc} and ~\ref{SSeta}.

{\bf Definition.} For $\eta\in Y^+$ we define the {\em extended local
convolution
diagram} $\CG\CS_\eta$ as the ind-scheme formed by the collections
$(\fV_\lambda,v_\lambda)_{\lambda\in X^+}$ of lattices in $V_\lambda\otimes\CK$
and vectors in $V_\lambda\otimes\CK$ such that

a) $(\fV_\lambda)_{\lambda\in X^+}\in\oCG_\eta$, or in other words,
$(\fV_\lambda)_{\lambda\in X^+}$ satisfies the condition ~\ref{loc}a)-c);

b) $v_\lambda\in\fV_\lambda$;

c) $(v_\lambda)_{\lambda\in X^+}\in\tCS$ (see ~\ref{SSeta}).

\subsection{Definition}
a) We define the projection $\bp:\ \CG\CS_\eta\lra\oCG_\eta$ as
$\bp(\fV_\lambda,v_\lambda)=(\fV_\lambda)$;

b) We define the map $\bq:\ \CG\CS_\eta\lra\CS$ as follows (cf. ~\ref{loc}a):
$$\bq(\fV_\lambda,v_\lambda)_{\lambda\in X^+}=
(z^{\langle\eta,\lambda\rangle}v_\lambda\in
z^{\langle\eta,\lambda\rangle}\fV_\lambda
\subset V_\lambda\otimes\CO)_{\lambda\in X^+}$$

\subsection{}
\label{lcd}
The torus $\bH_a$ acts in a natural way on $\CG\CS_\eta:\
h(\fV_\lambda,v_\lambda)=(\fV_\lambda,\lambda(h)v_\lambda)$. The action is
evidently free, and we denote the quotient by $\CG\bQ_\eta$,
the {\em local convolution diagram}. The map $\bp:\ \CG\CS_\eta\lra\oCG_\eta$
commutes with the action of $\bH_a$ (trivial on $\oCG_\eta$), so it descends
to the same named map $\bp:\ \CG\bQ_\eta\lra\oCG_\eta$.
The map $\bq:\ \CG\CS_\eta\lra\CS$ commutes with the action of $\bH_a$
(for the action on $\CS$ see ~\ref{SS}), so it descends to the same named
map $\bq:\ \CG\bQ_\eta\lra\bQ$.

The proalgebraic group $\bG(\CO)$ acts on $\CG\bQ_\eta$ and on $\bQ$, and
the map $\bq$ is equivariant with respect to this action.

\subsection{}
\label{cartes}
Let us compare the local convolution diagram with the global one.
For $\alpha\in Y$, taking formal expansion at $0\in C$ as in ~\ref{m},
we obtain the closed embedding $\fs:\ \CG\CQ^\alpha_\eta\hookrightarrow
\CG\bQ_\eta$. It is easy to see that the following diagram is cartesian:
$$
\begin{CD}
\CG\CQ^\alpha_\eta       @>\fs>>   \CG\bQ_\eta         \\
@V{\bq}VV                          @V{\bq}VV               \\
\CQ^{\eta+\alpha}    @>{\fs}>>             \bQ
\end{CD}
$$

\subsection{}
Recall the locally closed embedding $\CZ^\alpha\subset\CQ^\alpha$.

{\bf Definition.} We define the {\em restricted convolution diagram}
$\CG\CZ_\eta^\alpha\subset\CG\CQ_\eta^\alpha$ as the preimage
$\bq^{-1}(\CZ^{\eta+\alpha})$. The open subvariety $\CG\CZ_\eta^\alpha\cap
\oGQ_\eta^\alpha$
will be denoted by $\oGZ_\eta^\alpha$.
We will preserve
the notations $\bp,\bq$ (resp. $\obp,\bj$) for the restrictions of
these morphisms to $\CG\CZ_\eta^\alpha$ (resp. $\oGZ_\eta^\alpha$).

\subsection{}
\label{idiot}
We construct the {\em Beilinson-Drinfeld} avatar
$\CG\bZ^\alpha_\eta$ of the restricted convolution diagram $\CG\CZ_\eta^\alpha$.

{\bf Definition.} $\CG\bZ^\alpha_\eta$ is the space of triples
$(D,(\CU_\lambda)_{\lambda\in X^+},(\fU_\lambda)_{\lambda\in X^+})$
where $D\in\BA^{\eta+\alpha}$ is an effective $I$-colored divisor, and
$(\CU_\lambda)_{\lambda\in X^+},(\fU_\lambda)_{\lambda\in X^+}$ are the
collections of vector bundles on $\BP^1$ such that

a) $(\CU_\lambda)_{\lambda\in X^+}\in\oCG_\eta$, or in other words,
$(\CU_\lambda)_{\lambda\in X^+}$ satisfies the conditions ~\ref{closure} a)-c);

b) $(D,(\fU_\lambda)_{\lambda\in X^+})\in\bZ^{\eta+\alpha}$, or in other words,
$(D,(\fU_\lambda)_{\lambda\in X^+})$ satisfies the conditions ~\ref{PBD} a)-d);

c) $\iota_\lambda(\CV_\lambda^{\bN_-}(-\langle D,\lambda\rangle))\subset
\CU_\lambda(-\langle\eta,\lambda\rangle0)$
(notations of ~\ref{?}).

\subsubsection{}
The identification $\CG\CZ_\eta^\alpha=\CG\bZ_\eta^\alpha$ easily follows from
~\ref{PBD}. Under this identification, for
$(D,(\CU_\lambda)_{\lambda\in X^+},(\fU_\lambda)_{\lambda\in X^+})
\in\CG\bZ_\eta^\alpha$, we have
$$\bp(D,(\CU_\lambda)_{\lambda\in X^+},(\fU_\lambda)_{\lambda\in X^+})=
(\CU_\lambda)_{\lambda\in X^+}$$ and
$$\bq(D,(\CU_\lambda)_{\lambda\in X^+},(\fU_\lambda)_{\lambda\in X^+})=
(D,(\fU_\lambda)_{\lambda\in X^+})$$

\subsection{} We will introduce the {\em fine} stratifications of
$\CG\CQ_\eta^\alpha$ and $\CG\CZ_\eta^\alpha$ following the section 7.

\subsubsection{Fine stratification of $\CG\CQ_\eta^\alpha$}
\label{fine GQ}
We have $$\CG\CQ_\eta^\alpha=
\bigsqcup\ooGQ_\chi^\gamma\times(\BP^1-0)^{\beta-\gamma}_\Gamma$$
Here the union is taken over dominant $\chi\leq\eta$ in $Y^+$, arbitrary
$\gamma\leq\beta\leq\alpha\in Y$, and partitions $\Gamma\in\fP(\beta-\gamma)$.
Furthermore, $\ooGQ_\chi^\gamma\subset\oGQ_\chi^\gamma$ is an open subvariety
formed by all the collections $(\CU_\lambda,\fL_\lambda)_{\lambda\in X^+}$ such
that $(\CU_\lambda)\in\CG_\chi\subset\oCG_\chi$, and $\fL_\lambda$ is a
{\em line subbundle} in $\CU_\lambda$.

The stratum $\ooGQ_\chi^\gamma\times(\BP^1-0)^{\beta-\gamma}_\Gamma$ is formed
by the collections $(\CU_\lambda,\fL_\lambda)_{\lambda\in X^+}$ such that
$\bp(\CU_\lambda,\fL_\lambda)\in\CG_\chi\subset\oCG_\eta$; the normalization
(see ~\ref{strat M}) of $\fL$ in $\CU$ has degree $\gamma$; and the defect $D$
(see {\em loc. cit.}) of $\fL$ in $\CU$ equals $(\alpha-\beta)0+D'$ where
$D'\in(\BP^1-0)^{\beta-\gamma}_\Gamma$.

{\bf Proposition.} $\ooGQ^\alpha_\eta$ is smooth for arbitrary
$\alpha,\eta$, i.e. the fine strata are smooth.

{\em Proof} will be given in ~\ref{sometimes}.

\subsubsection{Fine Schubert stratification of $\CG\CQ_\eta^\alpha$}
\label{fineS GQ}
We have $$\CG\CQ_\eta^\alpha=
\bigsqcup\ooGQ_{w,\chi}^\gamma\times(\BP^1-0)^{\beta-\gamma}_\Gamma$$
Here the union is taken over dominant $\chi\leq\eta$ in $Y^+$,
representatives $w\in\CW_f/\CW_\chi$, arbitrary
$\gamma\leq\beta\leq\alpha\in Y$, and partitions $\Gamma\in\fP(\beta-\gamma)$.
Furthermore, $\ooGQ_{w,\chi}^\gamma\subset\oGQ_\chi^\gamma$
is an open subvariety
formed by all the collections $(\CU_\lambda,\fL_\lambda)_{\lambda\in X^+}$ such
that $(\CU_\lambda)\in\CG_{w,\chi}\subset\oCG_\chi$, and $\fL_\lambda$ is a
{\em line subbundle} in $\CU_\lambda$.

The stratum $\ooGQ_\chi^\gamma\times(\BP^1-0)^{\beta-\gamma}_\Gamma$ is formed
by the collections $(\CU_\lambda,\fL_\lambda)_{\lambda\in X^+}$ such that
$\bp(\CU_\lambda,\fL_\lambda)\in\CG_{w,\chi}\subset\oCG_\eta$;
the normalization of $\fL$ in $\CU$ has degree $\gamma$; and the defect $D$
(see {\em loc. cit.}) of $\fL$ in $\CU$ equals $(\alpha-\beta)0+D'$ where
$D'\in(\BP^1-0)^{\beta-\gamma}_\Gamma$.

{\bf Proposition.} $\ooGQ^\alpha_{w,\eta}$ is smooth for arbitrary
$\alpha,\eta,w$, i.e. the fine Schubert strata are smooth.

{\em Proof} will be given in ~\ref{sometimes}.

\subsubsection{Fine Schubert stratification of $\CG\CZ_\eta^\alpha$}
\label{fine GZ}
Similarly, we have $$\CG\CZ_\eta^\alpha=
\bigsqcup\oGZ_{w,\chi}^\gamma\times(\BC^*)^{\beta-\gamma}_\Gamma$$
Here the union is taken over dominant $\chi\leq\eta$ in $Y^+$,
representatives $w\in\CW_f/\CW_\chi$, arbitrary
$\gamma\leq\beta\leq\alpha\in Y$, and partitions $\Gamma\in\fP(\beta-\gamma)$.
Furthermore, $\oGZ_{w,\chi}^\gamma\subset\oGZ_\chi^\gamma$ is an open subvariety
formed by all the collections $(\CU_\lambda,\fL_\lambda)_{\lambda\in X^+}$ such
that $(\CU_\lambda)\in\CG_{w,\chi}\subset\oCG_\chi$,
and $\fL_\lambda$ is a {\em line subbundle} in $\CU_\lambda$.

The stratum $\oGZ_{w,\chi}^\gamma\times(\BC^*)^{\beta-\gamma}_\Gamma$ is formed
by the collections $(\CU_\lambda,\fL_\lambda)_{\lambda\in X^+}$ such that
$\bp(\CU_\lambda,\fL_\lambda)\in\CG_{w,\chi}\subset\oCG_\eta$; the normalization
of $\fL$ in $\CU$ has degree $\gamma$; and the defect $D$
of $\fL$ in $\CU$ equals $(\alpha-\beta)0+D'$ where
$D'\in(\BC^*)^{\beta-\gamma}_\Gamma$.

\subsubsection{}
\label{unwilling}
The reader unwilling to believe that $\oGZ^\gamma_{w,\chi}$
is smooth for arbitrary $\gamma\in Y$ may repeat the trick of
~\ref{awful}. Namely, one can replace $\CG\CZ^\alpha_\eta$ with an open
subvariety ${\ddot{\CG\CZ}}^\alpha_\eta$ formed by the union of the above
strata for sufficiently dominant $\gamma$; they are easily seen to be smooth.
Moreover, $\bq({\ddot{\CG\CZ}}^\alpha_\eta)\supset\ddZ^{\eta+\alpha}$, and
${\ddot{\CG\CZ}}^\alpha_\eta\supset\bq^{-1}(\ddZ^{\eta+\alpha})$.

\subsection{}
\label{fine local}
Let us introduce the {\em fine Schubert stratification} of the local
convolution diagram (see ~\ref{lcd}). Iwahori subgroup
$\bI\subset\bG(\CO)$ acts on $\tbQ$, and defines the {\em fine Schubert
stratification} of $\tbQ$ by Iwahori orbits $\dbQ_w^\alpha\subset\dbQ^\alpha,\
\alpha\in Y,w\in\CW_f$. Furthermore, we have
$$\CG\bQ_\eta=\bigsqcup\CG\dbQ_{w,\chi}^{-\alpha}$$
Here the union is taken over dominant $\chi\leq\eta$ in $Y^+$, representatives
$w\in\CW_f/\CW_\chi$, and $\alpha\in\BN[I]$. The stratum
$\CG\dbQ_{w,\chi}^{-\alpha}$ consists of collections
$(\fV_\lambda,v_\lambda)_{\lambda\in X^+}$ (vectors $v_\lambda\in\fV_\lambda$
are defined up to multiplication by $\BC^*$) such that

a) $(\fV_\lambda)_{\lambda\in X^+}\in\CG_{w,\chi}\subset\oCG_\eta$;

b) $z^{-\langle\alpha,\lambda\rangle}v_\lambda\in\fV_\lambda$ for all $\lambda$,
but $z^{-\langle\alpha,\lambda\rangle-1}v_\lambda\not\in\fV_\lambda$ for some
$\lambda$.

\section{Convolution}

\subsection{}
\label{hell}
Let $\fA,\fB$ be smooth varieties, and let $\fp:\ \fA\lra\fB$ be a map.
Suppose $\fA$ (resp. $\fB$) is equipped with a stratification $\fS$ (resp.
$\fT$), and $\fp$ is stratified with respect to the stratifications.
Let $\fR$ be another stratification of $\fB$, transversal to $\fT$.
Let $\CB$ (resp. $\CA$) be a perverse sheaf on $\fB$ (resp. $\fA$) smooth
along $\fR$ (resp. $\fS$). Let $b=\dim\fB$.


{\em Lemma.} a) $\CA\otimes\fp^*\CB[-b]$ is perverse;

b) $\CA\otimes\fp^*\CB[-b]=\CA\stackrel{!}{\otimes}\fp^!\CB[b]$;

c) Let $\ol{R}$ (resp. $\ol{S}$) be the closure of a stratum $R$ in $\fR$
(resp. $S$ in $\fS$). Then $\IC(\ol{S})\otimes\fp^*\IC(\ol{R})[-b]=
\IC(\ol{S}\cap\fp^{-1}\ol{R})$;

d) Let $R$ be a stratum of stratification $\fR$. Then $\fp^{-1}R$ is smooth.

{\em Proof.} a,b)
Let $g:\ \fG\lra\fA\times\fB$ denote the closed embedding of
the graph of $\fp$. The perverse sheaf $\CA\boxtimes\CB$ on $\fA\times\fB$
is smooth along the product stratification $\fS\times\fR$. The transversality
of $\fR$ and $\fT$ implies that the embedding $g$ is noncharacteristic with
respect to $\CA\boxtimes\CB$ (see ~\cite{ks}, Definition 5.4.12).
Furthermore, by definition, $\CA\otimes\fp^*\CB=g^*(\CA\boxtimes\CB)$, and
$\CA\stackrel{!}{\otimes}\fp^!\CB=g^!(\CA\boxtimes\CB)$. Now a) is nothing
else than ~\cite{ks}, Corollary 10.3.16(iii), while b) is nothing else than
~\cite{ks}, Proposition 5.4.13(ii). $\Box$

c) We consider $\ol{S}\times\ol{R}$ as a subvariety of $\fA\times\fB$.
Since $g:\ \fG\lra\fA\times\fB$ is noncharacteristic with respect to
$\IC(\ol{S})\boxtimes\IC(\ol{R})=\IC(\ol{S}\times\ol{R})$, we conclude that
$g^*\IC(\ol{S}\times\ol{R})[-b]=\IC(\fG\cap(\ol{S}\times\ol{R}))$.
It remains to note that $\fG\cap(\ol{S}\times\ol{R})=\ol{S}\cap\fp^{-1}\ol{R}$,
and $g^*\IC(\ol{S}\times\ol{R})=\IC(\ol{S})\otimes\fp^*\IC(\ol{R})$. $\Box$

d) We will view $\fp^{-1}R$ as a subscheme of $\fA$ (scheme-theoretic fiber
over $R$). Let $a\in\fp^{-1}R\subset\fA$.
The Zariski tangent space $T_a(\fp^{-1}R)$ equals $d\fp_a^{-1}(T_bR)$
where $b=\fp(a)$, and $d\fp_a:\ T_a\fA\lra T_b\fB$
stands for the differential of $\fp$ at $a$.
Let $\CT$ be a stratum of $\fT$ containing $b=\fp(a)$. Then $T_b\CT$ is
contained in $d\fp_a(T_a\fA)$ since $\fp$ is stratified with respect to $\fT$.
Furthermore, $T_bR+T_b\CT=T_b\fB$ due to the transversality of $\CT$ and $R$.
Hence $T_bR+d\fp_a(T_a\fA)=T_b\fB$. Hence $\dim(d\fp_a^{-1}(T_bR))=
\dim\fA-\dim\fB+\dim R$. We conclude that the dimension of the Zariski tangent
space $T_a(\fp^{-1}R)$ is independent of $a\in\fp^{-1}R$, and thus $\fp^{-1}R$
is smooth. $\Box$

\subsection{}
\label{j}
Consider the following cartesian diagram:
$$
\begin{CD}
\CG\CQ^\alpha_\eta       @>\bi>>    \fQ^\alpha        \\
@V{\bp}VV               @V{\bp}VV               \\
\oCG_\eta    @>{\bi}>>     \fM
\end{CD}
$$

{\bf Proposition.}
a) For a $\bG[[z]]$-equivariant perverse sheaf $\CF$ on $\oCG_\eta$ the sheaf
$\IC(\fQ^\alpha)\otimes\bp^*\bi_*\CF[-\dim\ufM^\eta]$ is supported on
$\CG\CQ^\alpha_\eta$ and is perverse;

b) $\IC(\CG\CQ_\eta^\alpha)=\IC(\fQ^\alpha)\otimes\bp^*\bi_*\IC(\oCG_\eta)
[-\dim\ufM^\eta]$.

{\em Proof.} a) Let us restrict the right column of the above diagram to
the open subscheme $\fM^\eta\subset\fM$, and take its quotient by $\bG^\eta$.
We obtain the cartesian diagram
$$
\begin{CD}
\CG\CQ^\alpha_\eta       @>\bi>>    \ufQ^{\alpha,\eta}        \\
@V{\bp}VV               @V{\bp}VV               \\
\oCG_\eta    @>{\bi}>>     \ufM^\eta
\end{CD}
$$
of schemes of finite type. Here the rows are closed embeddings, and
$\ufM^\eta$ is smooth. The stratification $\fS$ of $\fQ^\alpha$ is invariant
under the action of $\bG^\eta$, and descents to the same named quotient
stratification of $\ufQ^{\alpha,\eta}$. Similarly, the stratification of $\fM$
by the isomorphism type of $\bG$-torsors descents to the stratification
$\fT$ of $\ufM^\eta$. The sheaf $\bi_*\CF$ on $\ufM^\eta$ is smooth along
the stratification $\fR$ transversal to $\fT$. We have to prove that
$\IC(\ufQ^{\alpha,\eta})\otimes\bp^*\bi_*\CF[-\dim\ufM^\eta]$ is perverse.
In order to apply the Lemma ~\ref{hell}a) we only have to find an embedding
$\ufQ^{\alpha,\eta}\stackrel{u}{\hookrightarrow}\fA$ into a smooth scheme
such that the map $\bp$ and stratification $\fS$ extend to $\fA$. Then we
apply the Lemma ~\ref{hell}a) to the sheaf $\CA=u_*\IC(\ufQ^{\alpha,\eta})$
on $\fA$.

We will construct $\fA$ as a projective bundle over $\ufM^\eta$. The points of
$\ufM^\eta$ are the $\bG$-torsors $\CT$ over $C$ trivialized in some
infinitesimal neighbourhood of $\infty\in C$.
The points of $\ufQ^{\alpha,\eta}$ are the
$\bG$-torsors $\CT$ over $C$ trivialized
in some infinitesimal neighbourhood of $\infty\in C$ along with collections
of invertible subsheaves $\fL_\lambda\subset\CV^\CT_\lambda, \lambda\in X^+$,
satisfying Pl\"ucker relations.

Now, if $\alpha=\sum_{i\in I}a_ii$ is dominant enough,
the fiber of $\fA$ over $\CT\in\ufM^\eta$
is $\prod_{i\in I}\BP(\Gamma(C,\CV^\CT_{\omega_i}(a_i)))$. The map $u$
sends $(\CT,(\fL_\lambda\subset\CV^\CT_\lambda)_{\lambda\in X^+})$ to
$(\fL_{\omega_i}(a_i))\in
\prod_{i\in I}\BP(\Gamma(C,\CV^\CT_{\omega_i}(a_i)))$.

If $\alpha$ is not dominant enough, we first embed $\ufQ^{\alpha,\eta}$
into $\ufQ^{\beta+\alpha,\eta}$ for dominant enough $\beta$ as follows:
$(\CT,(\fL_\lambda\subset\CV^\CT_\lambda)_{\lambda\in X^+})\mapsto
(\CT,(\fL_\lambda(-\langle\beta,\lambda\rangle0)\subset
\CV^\CT_\lambda)_{\lambda\in X^+})$. Then we compose with the above projective
embedding of $\ufQ^{\beta+\alpha,\eta}$.

This completes the proof of a). $\Box$

b) We apply the Lemma ~\ref{hell}c) to $\ol{R}=\oCG_\eta,\
\ol{S}=\ufQ^{\alpha,\eta}$. $\Box$

\subsection{Proposition}
\label{+}
a) For a perverse sheaf $\CF$ in $\CP(\oCG_\eta,\bI)$ the sheaf
$\IC(\fQ^\alpha)\otimes\bp^*\bi_*\CF[-\dim\ufM^\eta]$ is supported on
$\CG\CQ^\alpha_\eta$ and is perverse;

b) $\IC(\CG\CQ_{w,\eta}^\alpha)=\IC(\fQ^\alpha)\otimes\bp^*\bi_*
\IC(\oCG_{w,\eta})[-\dim\ufM^\eta]$.

{\em Proof.} The same as the proof of the Proposition ~\ref{j}; we need only
to find a refinement $\fW$ of the stratification
$\fM=\bigsqcup_{\eta\in Y^+}\fM_\eta$ which would be transversal to the
Iwahori orbits $\CG_{w,\chi}$ in $\CG$.
Now $\fM=\bigsqcup_{\eta\in Y^+}\fM_\eta$ is the stratification by the
orbits of proalgebraic group $\bG[[z^{-1}]]$ acting naturally on $\fM$.
The desired refinement is the stratification by the orbits of subgroup
$\bI_-\subset\bG[[z^{-1}]]$ formed by the formal loops $g(z)\in\bG[[z^{-1}]]$
such that $g(\infty)\in\bB_-$. $\Box$

\subsubsection{}
Recall the notations of ~\ref{some}.

{\bf Conjecture.} a) The map $\obp:\ \oGQ^\alpha_\eta\lra\oCG_\eta$ is smooth
onto its image;

b) Up to a shift,
$\IC(\fQ^\alpha)\otimes\bp^*\bi_*\CF[-\dim\ufM^\eta]=\bj_{!*}\obp^*\CF$
for any perverse sheaf $\CF$ in $\CP(\oCG_\eta,\bI)$.

\subsection{Proof of the Propositions ~\ref{fine GQ} and ~\ref{fineS GQ}}
\label{sometimes}
We apply the Lemma ~\ref{hell}d) to the following situation:
$\fA=\overset{\circ}\ufQ{}^{\alpha,\eta}\subset\ufQ^{\alpha,\eta},\
\fB=\ufM^\eta,\ R=\CG_{w,\eta},\ \fp=\obp$. The stratification $\fT$ is defined
as follows. Recall the stratification
$\fW$ of $\fM$ by $\bI_-$-orbits introduced
in the proof of ~\ref{+}. It is invariant under the action of $\bG^\eta$
and descends to the desired stratification $\fT$ of $\ufM^\eta$ transversal
to $R$.

It remains to note that
$\fA=\overset{\circ}\ufQ{}^{\alpha,\eta}=\obp^{-1}(\fM^\eta)/\bG^\eta$
is smooth being a quotient by the free group action of an open subscheme
$\obp^{-1}(\fM^\eta)$ of the smooth scheme $\ofQ^\alpha$. Thus the
assumptions of ~\ref{hell}d) are in force, and we conclude that
$\ooGQ_{w,\eta}^\alpha=\obp^{-1}(\CG_{w,\eta})$ is smooth.

The proof of smoothness of $\ooGQ_\eta^\alpha$ is absolutely similar. $\Box$

\subsection{}
\label{z}

We define the following locally closed subscheme
$\fZ^\alpha\subset\fQ^\alpha$. Its points are
the $\bG$-torsors $\CT$ over $C$ trivialized
in the formal neighbourhood of $\infty\in C$ along with collections
of invertible subsheaves $\fL_\lambda\subset\CV^\CT_\lambda, \lambda\in X^+$,
satisfying Pl\"ucker relations plus two more conditions:

a) in some neighbourhood of $\infty\in C$ the invertible subsheaves
$\fL_\lambda\subset\CV^\CT_\lambda$ are line subbundles. Thus they may be
viewed as a reduction of $\CT$ to $\bB\subset\bG$ in this neighbourhood.
Since $\CT$ is trivialized in the formal neighbourhood of $\infty\in C$,
we obtain a map from this neighbourhood to the flag manifold $\bX$.

b) The value of the above map at $\infty\in C$ equals $\bB_-\in\bX$.

We have the following cartesian diagram:
$$
\begin{CD}
\CG\CZ^\alpha_\eta       @>\bi>>    \fZ^\alpha        \\
@V{\bp}VV                           @V{\bp}VV               \\
\oCG_\eta                @>{\bi}>>     \fM
\end{CD}
$$

{\bf Proposition.}
For a perverse sheaf $\CF\in\CP(\oCG_\eta,\bI)$ the sheaf
$\IC(\fZ^\alpha)\otimes\bp^*\bi_*\CF[-\dim\ufM^\eta]$ is supported on
$\CG\CZ^\alpha_\eta$ and is perverse.

{\em Proof.} Similar to the proof of the Proposition ~\ref{j}. $\Box$

\subsubsection{Remark}
\label{none}
Let us denote the embedding of $\fZ^\alpha$ into $\fQ^\alpha$ by $s$.
One can easily check that $\IC(\fZ^\alpha)=s^*\IC(\fQ^\alpha)[-\dim\bX]$.

%
%


\subsection{}
\label{isomor}
Recall the notations of ~\ref{sig}.

{\bf Proposition.} a) $\bq:\ \CG\CQ_\eta^\alpha\lra\CQ^{\eta+\alpha}$
(resp. $\bq:\ \CG\CZ_\eta^\alpha\lra\CZ^{\eta+\alpha}$) is proper;

b) Restriction of $\bq:\ \CG\CQ_\eta^\alpha\lra\CQ^{\eta+\alpha}$
(resp. $\bq:\ \CG\CZ_\eta^\alpha\lra\CZ^{\eta+\alpha}$) to $\qp^{\eta+\alpha}
\subset\CQ^{\eta+\alpha}$ (resp. to $\zp^{\eta+\alpha}\subset\CZ^{\eta+\alpha}$)
is an isomorphism.

{\em Proof.} a) is evident.

b) It suffices to consider the case of
$\bq:\ \CG\CQ^\alpha_\eta\lra\CQ^{\eta+\alpha}$.
Let $(\fL_\lambda)_{\lambda\in X^+}\in\CQ^{\eta+\alpha}$.
Then $(\CU_\lambda,\fL_\lambda')_{\lambda\in X^+}\in\bq^{-1}
(\fL_\lambda)_{\lambda\in X^+}$ iff for any $\lambda\in X^+$ we have
$\fL_\lambda'=\fL_\lambda(\langle\eta,\lambda\rangle0)$, and
$\CU_\lambda\supset\CV_\lambda(-\langle\eta,\lambda\rangle0)+\fL_\lambda'$.



Consider $\fL_\lambda=\CV_\lambda^{\bN_-}(-\langle\eta+\alpha,
\lambda\rangle\infty)$,
so that $\varphi=(\fL_\lambda)_{\lambda\in X^+}\in\qp^{\eta+\alpha}$.
Then $(\CU_\lambda,\fL_\lambda')_{\lambda\in X^+}\in\bq^{-1}
(\fL_\lambda)_{\lambda\in X^+}$ iff for any $\lambda\in X^+$ we have
$\fL_\lambda'=\fL_\lambda(\langle\eta,\lambda\rangle0)=
\CV_\lambda^{\bN_-}(\langle\eta,\lambda\rangle0-\langle\eta+
\alpha,\lambda\rangle\infty)$, and
$\CU_\lambda\supset\CV_\lambda(-\langle\eta,\lambda\rangle0)+\fL_\lambda'=
\CV_\lambda(-\langle\eta,\lambda\rangle0)+
\CV_\lambda^{\bN_-}(\langle\eta,\lambda\rangle0)$.

In other words, $(\CU_\lambda)_{\lambda\in X^+}$ lies in the intersection of
$\oCG_\eta$ with the semiinfinite orbit $T_\eta$ (see ~\ref{Iwasawa}). This
intersection consists exactly of one point (see ~\cite{mv}). Thus
$\bq^{-1}(\varphi)$ consists of one point.

Recall the cartesian diagram ~\ref{cartes}.
We have $\fs(\qp^{\eta+\alpha})\subset\dbQ^0$ (notations of ~\ref{QQ}),
in particular, $\fs(\varphi)\in\dbQ^0$. Since the map $\bq:\ \CG\bQ_\eta\lra
\bQ$ is $\bG(\CO)$-equivariant, and its fiber over $\fs(\varphi)$ consists
of one point, we conclude that $\bq$ is isomorphism over the $\bG(\CO)$-orbit
$\dbQ^0$. Since $\qp^{\eta+\alpha}=\fs^{-1}(\dbQ^0)$, applying the cartesian
diagram ~\ref{cartes}, we deduce that $\bq$ is isomorphism over
$\qp^{\eta+\alpha}.\ \Box$



\subsection{Proof of the Proposition ~\ref{GM}}
\label{later}
We are interested in the stalk of $\IC(\fQ^\alpha)$ at a point
$(\CT,(\fL_\lambda)_{\lambda\in X^+})\in\ofQ^\gamma\times
C^{\beta-\gamma}_\Gamma\subset\fQ^\alpha$. Suppose that the isomorphism class
of $\bG$-torsor $\CT$ equals $\eta\in Y^+$, i.e. $\CT\in\fM_\eta$.
The stalk in question evidently does not depend on a choice of $\CT\in\fM_\eta$
and the defect $D\in C^{\beta-\gamma}_\Gamma$. In particular, we may (and will)
suppose that $\CT\in\bi(\CG_\eta)$, and $D\in (C-0)^{\beta-\gamma}_\Gamma$.
Then one can see easily that the stalk in question is isomorphic, up to a shift,
to the stalk of Goresky-MacPherson sheaf $\IC(\CG\CQ^\alpha_\eta)$ at the
point $(\CT,(\fL_\lambda)_{\lambda\in X^+})\in\CG\CQ^\alpha_\eta$.

On the other hand, according to the Proposition ~\ref{isomor} b), the latter
stalk is isomorphic to the stalk of $\IC(\CQ^{\eta+\alpha})$ at the point
$\bq(\CT,(\fL_\lambda)_{\lambda\in X^+})$. This point has the same defect $D$.
Applying ~\ref{GMQ} we complete the proof of the Proposition ~\ref{GM}. $\Box$

\subsection{}
\label{semismall}
Recall that a map $\pi:\ \CX\lra\CY$ is called {\em dimensionally
semismall} if the following condition holds: let $\CY_m$ be the set
of all points $y\in\CY$ such that $\dim(\pi^{-1}y)\geq m$, then for $m>0$
we have $\codim_\CY\CY_m\geq2m$.
Let us define $\CX_m=\pi^{-1}\CY_m$. Then we can formulate an equivalent
condition of semismallness as follows: for any $m\geq0$ we have
$\codim_\CX\CX_m\geq m$.

Suppose $\CX$ (resp. $\CY$) is equipped with a stratification $\fS$
(resp. $\fT$), and $\pi$ is stratified with respect to $\fS$ and $\fT$.
Then $\pi$ is called {\em stratified semismall} (see ~\cite{mv}) if
$\pi$ is proper, and the restriction $\pi|_S$ to any stratum in $\fS$ is
dimensionally semismall. In this case
$\pi_*=\pi_!$ takes perverse sheaves on $\CX$ smooth along $\fS$ to perverse
sheaves on $\CY$ smooth along $\fT$ (see {\em loc. cit.}).

\subsection{}
\label{q}
Recall the fine Schubert stratifications of $\CQ^{\eta+\alpha}$
(resp. $\CG\CQ_\eta^\alpha$) introduced in ~\ref{fineS Q}
(resp. ~\ref{fineS GQ}). The map
$\bq:\ \CG\CQ_\eta^\alpha\lra\CQ^{\eta+\alpha}$ is stratified with respect
to these stratifications.

{\bf Proposition.}
$\bq:\ \CG\CQ_\eta^\alpha\lra\CQ^{\eta+\alpha}$ is stratified semismall.

%

{\em Proof} will use a few Lemmas.

\subsubsection{Lemma}
\label{lemma1}
$\bq:\ \CG\CQ_\eta^\alpha\lra\CQ^{\eta+\alpha}$ is dimensionally semismall.

{\em Proof.}
Recall that we have $\CQ^{\eta+\alpha}=\sqcup_{\beta\leq
\eta+\alpha}\qp^\beta$. It is enough to prove that for
$(\fL_\lambda)_{\lambda\in X^+}\in\qp^\beta$ we have
$\dim\bq^{-1}(\fL_\lambda)_{\lambda\in X^+}\leq |\eta+\alpha-\beta|$.

Let us start with the case $\fL_\lambda=
\CV_\lambda^{\bN_-}(\langle\beta-\eta-\alpha,\lambda\rangle0-
\langle\beta,\lambda\rangle\infty)$. Then, like in the Proposition
~\ref{isomor}, we have
$\bq^{-1}(\fL_\lambda)_{\lambda\in X^+}=\oCG_\eta\cap
\overline{T}_{\beta-\alpha}$, and
according to ~\cite{mv}, we have $\dim(\oCG_\eta\cap
\overline{T}_{\beta-\alpha})\leq
|\eta+\alpha-\beta|$.

Now $\qp^\beta$ is stratified by the defect: $\qp^\beta=
\sqcup^{\beta\geq\gamma\geq0}_{\Gamma\in\fP(\beta-\gamma)}
\qc^\gamma\times (C-0)^{\beta-\gamma}_\Gamma$, and $\bq$ is evidently
stratified with respect to this stratification. The point
$(\CV_\lambda^{\bN_-}(\langle\beta-\eta-\alpha,\lambda\rangle0-\langle\beta,
\lambda\rangle\infty))_{\lambda\in X^+}$ lies in the smallest (closed) stratum
$\gamma=0,\Gamma=\{\{\beta\}\}$. Since the dimension of preimage is a
lower semicontinuous function on $\qp^\beta$, we conclude that for
any point $(\fL_\lambda)_{\lambda\in X^+}\in\qp^\beta$ we have
$\dim\bq^{-1}(\fL_\lambda)_{\lambda\in X^+}\leq|\eta+\alpha-\beta|$.
$\Box$

\subsubsection{}
\label{lemma2}
Recall the fine stratification of $\CG\CQ^\alpha_\eta$
(resp. of $\CQ^{\eta+\alpha}$) introduced in ~\ref{fine GQ}
(resp. in ~\ref{fineS Q}).

{\em Lemma.}
$\bq:\ \CG\CQ_\eta^\alpha\lra\CQ^{\eta+\alpha}$ is stratified semismall
with respect to fine stratifications.

{\em Proof.} We consider a fine stratum
$\ooGQ_\chi^\gamma\times(\BP^1-0)^{\beta-\gamma}_\Gamma\subset
\CG\CQ^\alpha_\eta$ (see ~\ref{fine GQ}).
Temporarily we will write $\bq^\alpha_\eta$ for $\bq$ to stress its dependence
on $\eta$ and $\alpha$.

The restriction of $\bq^\alpha_\eta$ to the stratum
$\ooGQ_\chi^\gamma\times(\BP^1-0)^{\beta-\gamma}_\Gamma$ decomposes into
the following composition of morphisms:
$$\ooGQ_\chi^\gamma\times(\BP^1-0)^{\beta-\gamma}_\Gamma
\stackrel{a\times\id}{\hookrightarrow}
\CG\CQ_\chi^\gamma\times(\BP^1-0)^{\beta-\gamma}_\Gamma
\stackrel{\bq^\gamma_\chi\times\id}{\lra}
\CQ^{\chi+\gamma}\times(\BP^1-0)^{\beta-\gamma}_\Gamma
\stackrel{b}{\hookrightarrow}
\CQ^{\eta+\alpha}$$
Here $a$ is the open inclusion; and $b(\fL,D')=
\fL((\beta-\alpha+\chi-\eta)0-D')$.

Now $\bq^\gamma_\chi\times\id$ is semismall according to the Lemma
~\ref{lemma1}. This completes the proof of the Lemma. $\Box$

\subsubsection{Lemma}
\label{lemma3}
The restriction of $\bq$ to the fine Schubert stratum
$\ooGQ^\alpha_{w,\eta}\subset\CG\CQ^\alpha_\eta$ is dimensionally
semismall for any $w\in\CW_f/\CW_\eta$.

{\em Proof.} Let $\bK\subset\bI\subset\bG(\CO)$
denote the first congruence subgroup
formed by the loops $g(z)\in\bG(\CO)$ such that $g(0)=1$.
The point $\eta(z)\in\CG_{e,\eta}\subset\CG_\eta$ was introduced in
~\ref{Iwahori}. For a positive integer $m$ the subset $(\CG\CQ^\alpha_\eta)_m$
(with respect to $\bq:\ \CG\CQ^\alpha_\eta\lra\CQ^{\eta+\alpha}$)
was introduced in ~\ref{semismall}.

{\em Claim 1.} $\codim_{\ooGQ_\eta^\alpha}(\ooGQ_\eta^\alpha)_m=
\codim_{\obp^{-1}(g\cdot\bK\cdot\eta(z))}[\obp^{-1}(g\cdot\bK\cdot\eta(z))
\cap(\ooGQ^\alpha_\eta)_m]$ for any $m\geq0$ and $g\in\bG$.

In effect, due to the $\bG$-equivariance of $\obp$ and $\bq$, the RHS does
not depend on a choice of $g\in\bG$, so it suffices to consider $g=e$.

The stabilizer of $\eta(z)$ in $\bG$ is nothing else than the parabolic
subgroup $\bP(I_\eta)$ introduced in ~\ref{Iwahori}.
We have
$$\CG_\eta=\bG\times_{\bP(I_\eta)}[\bK\cdot\eta(z)];\
\ooGQ_\eta^\alpha=\bG\times_{\bP(I_\eta)}[\obp^{-1}(\bK\cdot\eta(z))];\
(\ooGQ_\eta^\alpha)_m=\bG\times_{\bP(I_\eta)}
[\obp^{-1}(\bK\cdot\eta(z))\cap(\ooGQ^\alpha_\eta)_m]$$
The Claim follows.

{\em Claim 2.}
$\codim_{\obp^{-1}(g\cdot\bK\cdot\eta(z))}
[\obp^{-1}(g\cdot\bK\cdot\eta(z))\cap(\ooGQ^\alpha_\eta)_m]=
\codim_{\ooGQ^\alpha_{w,\eta}}
[\ooGQ^\alpha_{w,\eta}\cap(\ooGQ^\alpha_\eta)_m]$
for any $m\geq0,\ g\in\bG$, and $w\in\CW_f/\CW_\eta$.

In effect, let us choose $g$ in the normalizer of $\bH$ representing $w$.
Let us denote by $\bP_w$ the intersection $\bP(I_\eta)\cap g\bP(I_\eta)g^{-1}$.
Then we have
$$\ooGQ^\alpha_{w,\eta}=\bP(I_\eta)\times_{\bP_w}
[\obp^{-1}(g\cdot\bK\cdot\eta(z))];\
\ooGQ^\alpha_{w,\eta}\cap(\ooGQ^\alpha_\eta)_m=\bP(I_\eta)\times_{\bP_w}
[\obp^{-1}(g\cdot\bK\cdot\eta(z))\cap(\ooGQ^\alpha_\eta)_m]$$
The Claim follows.

Comparing the two Claims we obtain
$$\codim_{\ooGQ^\alpha_{w,\eta}}(\ooGQ^\alpha_{w,\eta})_m\geq
\codim_{\ooGQ^\alpha_{w,\eta}}[\ooGQ^\alpha_{w,\eta}\cap(\ooGQ^\alpha_\eta)_m]=
\codim_{\ooGQ^\alpha_\eta}(\ooGQ^\alpha_\eta)_m\geq m$$
The last inequality holds by the virtue of the Lemma ~\ref{lemma1}.
This completes the proof of the Lemma. $\Box$

\subsubsection{} Now we are ready to finish the proof of the Proposition.
It remains to show that the restriction of $\bq$ to any fine Schubert stratum
is dimensionally semismall. It follows from the Lemma ~\ref{lemma3}
in the same way as the Lemma ~\ref{lemma2} followed from the
Lemma ~\ref{lemma1} (``twisting by defect'').
This completes the proof of the Proposition. $\Box$

\subsection{Corollary}
\label{exact}
The functor $\bq_*=\bq_!$ takes perverse sheaves on $\CG\CQ_\eta^\alpha$
smooth along the fine Schubert stratification to perverse sheaves on
$\CQ^{\eta+\alpha}$ smooth along the fine Schubert stratification. $\Box$

\subsection{} A few remarks are in order.

\subsubsection{Remark}
Recall the fine Schubert stratification of the local convolution diagram
introduced in ~\ref{lcd} and ~\ref{fine local}.
The arguments used in the proof of the Proposition ~\ref{q} along with the
Lemma ~\ref{codime} show that the map $\bq:\ \CG\bQ_\eta\lra\bQ$ is stratified
semismall with respect to the fine Schubert stratification.

\subsubsection{Remark}
The same arguments as in the proof of Lemma ~\ref{lemma3} show that the
convolution $\CA*\CB$ of perverse sheaves on $\CG$ (see
~\cite{lus}, or ~\cite{g}, ~\cite{mv})
is perverse if $\CB$ is $\bG(\CO)$-equivariant.
In the particular case $\CA\in\CP(\CG,\bI)$ this was also proved by
G.Lusztig in ~\cite{l5} using calculations in the affine Hecke algebra.

\subsection{Theorem}
\label{tough}
Let $\eta\in Y^+,\alpha\in\BN[I]$.
Consider the following diagram:
$$
\begin{CD}
\oCG_\eta  @<{\bp}<< \CG\CQ^\alpha_\eta    @>{\bq}>> \CQ^{\eta+\alpha}       \\
@V{\bi}V{\subset}V       @V{\bi}V{\subset}V               @.            \\
\fM         @<\bp<<        \fQ^\alpha                  @.         {}
\end{CD}
$$

For a perverse sheaf $\CF$ in $\CP(\oCG_\eta,\bI)$, the sheaf
$\bc^\alpha_\CQ(\CF):=\bq_*(\IC(\fQ^\alpha)\otimes\bp^*\CF)[-\dim\ufM^\eta]$ on
$\CQ^{\eta+\alpha}$ is perverse and smooth along
the fine Schubert stratification.

{\em Proof.} By the virtue of ~\ref{+} we know that
$\IC(\fQ^\alpha)\otimes\bp^*\CF[-\dim\ufM^\eta]$ is a perverse sheaf on
$\CG\CQ^\alpha_\eta$. In order to apply the Corollary ~\ref{exact} we have
to check that
$\IC(\fQ^\alpha)\otimes\bp^*\CF[-\dim\ufM^\eta]$ is smooth along the fine
Schubert stratification. The sheaf $\bp^*\CF$ is evidently smooth along
the fine Schubert stratification. The sheaf $\IC(\fQ^\alpha)$ is constant
along the stratification by defect (see ~\ref{GM}), hence
$\bi^*\IC(\fQ^\alpha)$ is smooth along the fine Schubert stratification.
This completes the proof of the Theorem. $\Box$

\subsection{Conjecture}
Let $\eta\in Y^+,\alpha\in Y$ be such that $\eta+\alpha\in\BN[I]$.
Consider the following diagram:
$$
\begin{CD}
\oCG_\eta  @<{\bp}<< \CG\CQ^\alpha_\eta    @>{\bq}>> \CQ^{\eta+\alpha}       \\
@V{\bi}V{\subset}V       @V{\bi}V{\subset}V               @.            \\
\fM         @<\bp<<        \fQ^\alpha                  @.         {}
\end{CD}
$$

For a perverse sheaf $\CF$ in $\CP(\oCG_\eta,\bI)$, the sheaf
$\bq_*(\IC(\fQ^\alpha)\otimes\bp^*\CF)[-\dim\ufM^\eta]$ on
$\CQ^{\eta+\alpha}$ is perverse and smooth along
the fine Schubert stratification.

\subsection{Corollary}
\label{bunk}
Let $\eta\in Y^+,\alpha\gg0$.
Consider the following diagram:
$$
\begin{CD}
\oCG_\eta  @<{\bp}<< {\ddot{\CG\CZ}}^\alpha_\eta @>\bq>> \ddZ^{\eta+\alpha} \\
@V{\bi}V{\subset}V       @V{\bi}V{\subset}V               @.            \\
\fM         @<\bp<<     \fZ^\alpha                  @.         {}
\end{CD}
$$
(notations of ~\ref{unwilling} and ~\ref{z}).

For a perverse sheaf $\CF\in\CP(\oCG_\eta,\bI)$, the sheaf
$\bc^\alpha_\CZ(\CF):=\bq_*(\IC(\fZ^\alpha)\otimes\bp^*\CF)[-\dim\ufM^\eta]$
on $\ddZ^{\eta+\alpha}$ is perverse and smooth along
the fine Schubert stratification.

{\em Proof.} Let us denote by $s$ the locally closed embedding
$\ddZ^{\eta+\alpha}\stackrel{s}{\hookrightarrow}\CQ^{\eta+\alpha}$.
Also, temporarily, let us denote the maps $\bp$ and $\bq$ from the diagram
~\ref{tough} (resp. ~\ref{bunk}) by $\bp_\CQ$ and $\bq^\CQ$ (resp. $\bp_\CZ$
and $\bq^\CZ$) to stress their difference.

Then we have
$\bq^\CZ_*(\IC(\fZ^\alpha)\otimes\bp_\CZ^*\CF)[-\dim\ufM^\eta]=
s^*\bq^\CQ_*(\IC(\fQ^\alpha)\otimes\bp_\CQ^*\CF)[-\dim\ufM^\eta-\dim\bX]$
(cf. ~\ref{none}).
We also have $\codim_{\CQ^{\eta+\alpha}}\ddZ^{\eta+\alpha}=\dim\bX$, and
the fine Schubert strata in $\ddZ^{\eta+\alpha}$ are intersections of
$\ddZ^{\eta+\alpha}$ with the fine Schubert strata in $\CQ^{\eta+\alpha}$.
One can check readily that the functor $s^*[-\dim\bX]$ takes perverse
sheaves on $\CQ^{\eta+\alpha}$ smooth along the fine Schubert stratification
to perverse sheaves on $\ddZ^{\eta+\alpha}$ smooth along the fine Schubert
stratification. The application of ~\ref{tough} completes the proof of the
Corollary. $\Box$

\subsection{Conjecture}
Let $\eta\in Y^+,\alpha\in Y$ be such that $\eta+\alpha\in\BN[I]$.
Consider the following diagram:
$$
\begin{CD}
\oCG_\eta  @<{\bp}<< \CG\CZ^\alpha_\eta    @>{\bq}>> \CZ^{\eta+\alpha}       \\
@V{\bi}V{\subset}V       @V{\bi}V{\subset}V               @.            \\
\fM         @<\bp<<     \fZ^\alpha                  @.         {}
\end{CD}
$$

For a perverse sheaf $\CF$ in $\CP(\oCG_\eta,\bI)$, the sheaf
$\bq_*(\IC(\fZ^\alpha)\otimes\bp^*\CF)[-\dim\ufM^\eta]$ on
$\CZ^{\eta+\alpha}$ is perverse and smooth along
the fine Schubert stratification.


\subsection{}

Now we will compare $\bc_\CZ^\alpha(\CF)$ for a fixed
$\CF\in\CP(\oCG_\eta,\bI)$
and various $\alpha$. Recall the notations of ~\ref{snop} and ~\ref{awful}.

{\bf Proposition.} For any $\beta,\gamma\in\BN[I], \varepsilon>0,$ there is a
factorization isomorphism
$$\bc_\CZ^{\beta+\gamma-\eta}\CF|_{\ddZ^{\beta,\gamma}_{\Ue,\Upe}}\iso
\bc_\CZ^{\beta-\eta}\CF|_{\ddZ^\beta_\Ue}\boxtimes
\IC^\gamma|_{\ddZ^\gamma_\Upe}$$

{\em Proof.} Follows easily from ~\ref{idiot}. $\Box$

\subsection{} The above Proposition shows that we can organize the collection
$(\bc_\CZ^{\alpha-\eta}\CF)$ for $\alpha\in\BN[I]$ into a snop $\bc_\CZ\CF$.
Namely, we set the support estimate $\chi(\bc_\CZ\CF)=\eta,\
(\bc_\CZ\CF)^\alpha_\eta=\bc_\CZ^{\alpha-\eta}\CF$. This way we obtain an
exact functor $\bc_\CZ:\ \CP(\CG,\bI)\lra\PS$.


\section{Examples of convolution}

\subsection{}
\label{stalks}
Let $\CF$ be a perverse sheaf in $\CP(\oCG_\eta,\bI)$.
For $\chi\leq\eta, w\in\CW_f/\CW_\chi$ we have $\CG_{w,\chi}\subset\oCG_\eta$.
The sheaf $\CF$ is constant along $\CG_{w,\chi}$, and we denote by
$\CF_{w,\chi}$ its stalk at any point in $\CG_{w,\chi}$.

{\em Lemma.} The stalk of $\IC(\fQ^\alpha)\otimes\bp^*\CF$ at any
point in a fine Schubert stratum $\ooGQ^\gamma_{w,\chi}\times
(\BC^*)^{\beta-\gamma}_\Gamma\subset\CG\CQ_\eta^\alpha$ (see ~\ref{fineS GQ})
equals $\CF_{w,\chi}\otimes\IC^{\alpha-\beta}_\Gamma$ (see ~\ref{GMZ}).

{\em Proof.} Follows immediately from the Proposition ~\ref{GM}.
$\Box$

\subsection{}
\label{Satake}
Let $\bG^L$ be the Langlands dual group. Its character lattice
coincides with $Y$, and the dominant characters are exactly $Y^+$. For
$\eta\in Y^+$ we denote by $W_\eta$ the irreducible $\bG^L$-module with the
highest weight $\eta$. For $\chi\in Y$ we denote by $_{(\chi)}W_\eta$ the
weight $\chi$-subspace of $W_\eta$.

Let $\IC(\oCG_\eta)$ denote the Goresky-MacPherson sheaf of $\oCG_\eta$.
A natural isomorphism $H^\bullet(\oCG_\eta,\IC(\oCG_\eta))\iso W_\eta$ is
constructed in ~\cite{mv}.

Recall that for $\chi\in Y,\ w\in\CW_f,$ the irreducible snop $\CL(w,\chi)$
was introduced in ~\ref{CL}. The following result was suggested by V.Ginzburg.

{\bf Theorem.} There is a natural isomorphism of snops:
$$\bc_\CZ\IC(\oCG_\eta)\iso\bigoplus_{\chi\in Y}\ _{(\chi)}W_\eta\otimes
\CL(w_0,\chi)$$

{\em Proof.} It is a reformulation of the main result of ~\cite{mv}.
In effect, by the Proposition ~\ref{j}b) we know that
$\IC(\fQ^\alpha)\otimes\bp^*\IC(\oCG_\eta)[-\dim\ufM^\eta]=
\IC(\CG\CQ_\eta^\alpha)$. So
we have to prove that $\bq_*\IC(\CG\CQ_\eta^\alpha)=
\bigoplus_{0\leq\beta\leq\eta+\alpha}\
_{(\beta-\alpha)}W_\eta\otimes\IC(\CQ^\beta)$. Here we make use of the
filtration $\CQ^{\eta+\alpha}=\bigcup_{0\leq\beta\leq\eta+\alpha}\CQ^\beta$
subject to the stratification $\CQ^{\eta+\alpha}=
\bigsqcup_{0\leq\beta\leq\eta+\alpha}\qp^\beta$ by the defect at $0\in C$.

We know that $\bq$ is proper, semismall, and stratified with respect to the
above stratification. By the Decomposition Theorem (see ~\cite{bbd}), we
have {\em a priori}
$\bq_*\IC(\CG\CQ_\eta^\alpha)=\bigoplus_{0\leq\beta\leq\eta+\alpha}
L_\beta\otimes\IC(\CQ^\beta)$ for some vector spaces $L_\beta$. To identify
$L_\beta$ with $_{(\beta-\alpha)}W_\eta$ it suffices to compute the stalks at
$\phi=(\fL_\lambda)_{\lambda\in X^+}\in\qp^\beta$ where $\fL_\lambda=
\CV_\lambda^{\bN_-}(\langle\beta-\eta-\alpha,\lambda\rangle\cdot0-
\langle D,\lambda\rangle)$ for some $D\in(\BP^1-0)^\beta_\Gamma$.

As in the proof of ~\ref{lemma1} we have
$\bq^{-1}(\phi)=\oCG_\eta\cap\overline{T}_{\beta-\alpha}=
\bigsqcup_{\gamma\geq0}
(\oCG_\eta\cap T_{\beta-\alpha+\gamma})$.
According to the Lemma ~\ref{stalks}, we have
$\IC(\CG\CQ_\eta^\alpha)|_{\oCG_\eta\cap T_{\beta-\alpha+\gamma}}=
\IC(\oCG_\eta)|_{\oCG_\eta\cap T_{\beta-\alpha+\gamma}}\otimes
\IC^\gamma_\Gamma$. According to ~\cite{mv}, we have
$H^\bullet_c(\oCG_\eta\cap T_{\beta-\alpha+\gamma},\IC(\oCG_\eta))=\
_{(\beta-\alpha+\gamma)}W_\eta$. Due to the parity vanishing (see
{\em loc. cit.} and ~\ref{parity}), the spectral sequence computing
$H^\bullet(\oCG_\eta\cap\overline{T}_{\beta-\alpha},\IC(\CG\CQ_\eta^\alpha))$
collapses and gives
$H^\bullet(\oCG_\eta\cap\overline{T}_{\beta-\alpha},\IC(\CG\CQ_\eta^\alpha))=
\bigoplus_{0\leq\gamma\leq\eta+\alpha-\beta}\
_{(\beta-\alpha+\gamma)}W_\eta\otimes\IC^0_{\{\{\gamma\}\}}=
\bigoplus_{0\leq\gamma\leq\eta+\alpha-\beta}\
_{(\beta-\alpha+\gamma)}W_\eta\otimes\IC(\CQ^{\beta+\gamma})_\phi$.

This completes the proof of the Theorem. $\Box$

\end{document}